\documentclass{article}%
\usepackage{amssymb}
\usepackage{graphicx}
\usepackage{ifpdf}
\usepackage{srcltx}
\usepackage{url}
\usepackage{hyperref}
\usepackage{cite}
\usepackage{color}
\usepackage{multirow}
\usepackage{caption}
\usepackage{epstopdf}
\usepackage{letltxmacro}
\usepackage[cmex10]{amsmath}
\usepackage{amssymb}
\usepackage{algorithmic}
\usepackage{array}
\usepackage{eqparbox}
\usepackage[tight,footnotesize]{subfigure}
\usepackage{tabularx}
\usepackage{hyphenat}
\usepackage{amsfonts}
\usepackage{amsmath}
\usepackage[left=0.6in, right=0.6in, bottom=1.1in,top=1.1in]{geometry}%
\usepackage{graphicx}
\providecommand{\U}[1]{\protect\rule{.1in}{.1in}}

\begin{document}

\author{M. S. Eliwa$^{1}$\ and M. El-Morshedy$^{1,}$\thanks{Corresponding author:
mah\_elmorshedy@mans.edu.eg}\\
$^{1}$Mathematics Department, Faculty of Science, Mansoura \\
University,\ Mansoura, Pin 35516, Egypt.}
\title{Bivariate Discrete Inverse Weibull Distribution}
\date{}
\maketitle

\begin{abstract}
In this paper, we propose a new class of bivariate distributions, called the
bivariate discrete inverse Weibull (BDsIW) distribution, whose marginals are
discrete inverse Weibull (DsIW) distributions. Some statistical and
mathematical properties are presented. The maximum likelihood method is used
for estimating the model parameters. Simulations are presented to verify the
performance of the direct maximum likelihood estimation. Finally, two real
data sets are analyzed for illustrative purposes.



\end{abstract}

\section{Introduction}

The Weibull (W) distribution is one of the most popular and widely used
distributions for failure time in reliability theory (see, Weibull (1951)).
The cumulative distribution function (CDF) of W distribution is given by \ \
\begin{equation}
\Pi(x;\nu,\zeta)=1-e^{-\nu x^{\zeta}};\ \ x>0,
\end{equation}
where $\nu>0$ is scale parameter and $\zeta>0$ is shape parameter. Clearly,
the exponential (E) distribution and the Rayleigh (R) distribution are special
cases for $\zeta=1$ and $\zeta=2$ respectively. Unfortunately, the shape of
the hazard rate function (HRF) of W distribution can be only increasing,
decreasing or constant. So, more modifications and generalizations of W
distribution are presented in the statistical literature to describe various
phenomena in different fields, because in many applications, empirical hazard
rate curves often exhibit non-monotonic shapes such as a bathtub, unimodal and
others. For example:

\begin{enumerate}
\item Keller et al. (1985) proposed inverse Weibull (IW) distribution. The CDF
of IW distribution is given by%
\begin{equation}
\Pi(x;\nu,\zeta)=e^{-\nu x^{-\zeta}};\ \ x>0.
\end{equation}

\item Lai et al. (2003) introduced modified Weibull (MW) distribution. The CDF
of MW distribution is given by%
\begin{equation}
\Pi(x;\nu,\zeta,\lambda)=1-e^{-\nu x^{\zeta}e^{\lambda x}};\ \ x>0,
\end{equation}
where $\lambda>0$ is an accelerating parameter. The exponentiated MW
distribution was proposed by Jalmar et al. (2008).

\item Bebbington et al. (2007) proposed flexible Weibull (FxW) distribution.
The CDF of FxW distribution is given by%
\begin{equation}
\Pi(x;\nu,\gamma)=1-e^{-e^{\nu x-\frac{\gamma}{x}}};\ \ x>0,
\end{equation}
where $\gamma>0$ is scale parameter. The exponentiated FxW distribution was
presented by El-Gohary et al. (2015a), the inverse flexible Weibull (IFxW)
distribution was proposed by El-Gohary et al. (2015b) and the exponentiated of
it was studied by El-Morshedy et al. (2017).

\item Cordeiro et al. (2013) introduced exponential-Weibull (E-W)
distribution. The CDF of E-W distribution is given by%
\begin{equation}
\Pi(x;\gamma,\nu,\beta)=1-e^{-\gamma x-\nu x^{\beta}};\ \ x>0,
\end{equation}
where $\beta\in(0,\infty)-\{1\}$ is shape parameter.

\item Nadarajah et al. (2013) proposed exponentiated Weibull (EW)
distribution. The CDF of EW distribution is given by%
\begin{equation}
\Pi(x;\nu,\zeta,\eta)=\left(  1-e^{-\left(  \nu x\right)  ^{\zeta}}\right)
^{\eta};\ \ x>0,
\end{equation}
where $\eta>0$ is shape parameter.

\item El-Bassiouny et al. (2017) introduced exponentiated generalized
Weibull-Gompertz (EGW-Gz) distribution. The CDF of EGW-Gz distribution is
given by%
\begin{equation}
\Pi(x;\nu,\zeta,\lambda,\eta,\rho)=\left(  1-e^{-\nu x^{\zeta}\left(
e^{\lambda x^{\eta}}-1\right)  }\right)  ^{\rho};\ \ x>0,
\end{equation}
where $\rho>0$ is shape parameter. The mixure of 2-EGW-Gz distribution was
studied by El-Bassiouny et al. (2016).
\end{enumerate}

Moreover, some discrete versions of E, R, W distributions and its
generalizations have been presented in the literature because in several
cases, lifetimes need to be recorded on a discrete scale rather than on a
continuous analogue. So, discretizing continuous distributions has received
much attention in the literature. For example:

\begin{enumerate}
\item Toshio and Shunji (1975) introduced discrete Weibull (DsW) distribution.
The probability mass function (PMF) of DsW distribution is given by \
\begin{equation}
\pi(x;\theta,\zeta)=\theta^{x^{\zeta}}-\theta^{\left(  x+1\right)  ^{\zeta}%
};\ \ x\in\mathbf{%
\mathbb{N}
}_{0}=\left\{  0,1,2,3,...\right\}  ,
\end{equation}
where $0<\theta<1$. Clearly, the discrete Rayleigh (DsR) distribution is a
special case for $\zeta=2,$ which was presented by Dilip (2004).

\item G\'{o}mez-D\'{e}niz (2010) proposed generalization of geometric (GGo)
distribution. The PMF of GGo distribution is given by
\begin{equation}
\pi(x;\theta,\gamma)=\frac{\gamma\theta^{x}\left(  1-\theta\right)  }{\left(
1-\left[  1-\gamma\right]  \theta^{x+1}\right)  \left(  1-\left[
1-\gamma\right]  \theta^{x}\right)  };\ \ x\in\mathbf{%
\mathbb{N}
}_{0}.
\end{equation}
The GGo distribution reduces to geometric or discrete exponential (DsE)
distribution when $\gamma=1.$

\item Jazi et al. (2010) introduced DsIW distribution. The PMF of DsIW
distribution is given by%
\begin{equation}
\pi(x;\theta,\zeta)=\theta^{\left(  x+1\right)  ^{-\zeta}}-\theta^{x^{-\zeta}%
};\ \ x\in\mathbf{%
\mathbb{N}
}_{0}.
\end{equation}

\item Vahid et al. (2013) proposed discrete generalized exponential type II
(DsGE-T2) distribution. The PMF of DsGE-T2 distribution is given by%
\begin{equation}
\pi(x;\theta,\zeta)=\left(  1-\theta^{x+1}\right)  ^{\zeta}-\left(
1-\theta^{x}\right)  ^{\zeta};\ \ x\in\mathbf{%
\mathbb{N}
}_{0}.
\end{equation}

\item Vahid and Hamid (2015a) introduced exponentiated discrete Weibull (EDsW)
distribution. The PMF of EDsW distribution is given by
\begin{equation}
\pi(x;\theta,\rho,\zeta)=\left(  1-\theta^{\left(  x+1\right)  ^{\rho}%
}\right)  ^{\zeta}-\left(  1-\theta^{x^{\rho}}\right)  ^{\zeta};\ \ x\in
\mathbf{%
\mathbb{N}
}_{0}.
\end{equation}

\item Vahid et al. (2015b) proposed discrete beta exponential (DsBE)
distribution. The PMF of DsBE distribution is given by%
\begin{equation}
\pi(x;\theta,\gamma,\zeta)=c\theta^{\gamma\left(  x-1\right)  }\left(
1-\theta^{x}\right)  ^{\zeta-1};\ \ x\in\mathbf{%
\mathbb{N}
}_{0}-\left\{  0\right\}  ,
\end{equation}
where $c^{-1}=\sum_{j=0}^{\infty}\frac{\left(  -\theta\right)  ^{j}}%
{1-\theta^{\gamma+j}}\frac{\left(  \zeta-1\right)  \left(  \zeta-2\right)
...\left(  \zeta-j\right)  }{j!}$.
\end{enumerate}

On the other hand, in many practical situations, it is important to consider
different bivariate continuous and discrete distributions that could be used
to model bivariate lifetime data in many fields. So, several bivariate
continuous and discrete distributions are available in the statistical
literature. For example, Lee (1997), Karlis and Ntzoufras (2000), Wu and Yuen
(2003), Yuen et al. (2006), Sarhan and Balakrishnan (2007), Kundu and Dey
(2009), Morata (2009), Kundu and Gupta (2009), Ong and Ng (2013), Balakrishnan
and Shiji (2014), Lee and Cha (2015), Rasool and Akbar (2016), Hiba (2016),
El-Bassiouny et al. (2016), El-Gohary et al. (2016), Vahid and Kundu (2017),
Mohamed et al. (2017), Kundu and Vahid (2018), El-Morshedy and Khalil (2018)
among others. An excellent encyclopedic survey of various continuous and
discrete bivariate distributions can be found in Balakrishnan and Lai (2009)
and Johnson et al. (1997) respectively.

In this regard, we focus the aim of this paper on presenting a flexible
discrete bivariate distribution called BDsIW distribution, which can be
usefully applied not only by statisticians, but also by data analysis in many
different disciplines, such as sports, engineering, and medical applications.
The proposed discrete model can be obtained from 3-independent DsIW
distributions by using the maximization method as suggested by Lee and Cha
(2015). The main reasons for introducing BDsIW distribution are:

\begin{enumerate}
\item The proposed model is a very flexible bivariate discrete distribution,
and its joint PMF can take different shapes depending on the parameter values.

\item The generation from the proposed model is straight forward. So, the
simulation experiments can be performed quite conveniently.

\item The marginals of the proposed model are DsIW distributions. Hence, the
marginals are able to analyze the hazard rates in the discrete case.

\item The DsE and DsR distributions are special cases from the proposed model.
\end{enumerate}

\section{The BDsIW Distribution and Its Statistical Properties}

\subsection{Definition and interpretations}

Suppose $W_{1}\sim$ DsIW($\theta_{1},\zeta$), $W_{2}\sim$ DsIW($\theta
_{2},\zeta$) and $W_{3}\sim$ DsIW($\theta_{3},\zeta$) and they are
independently distributed. If $X_{d}=\max(W_{d},W_{3});d=1,2$, then we can say
that the bivariate vector $\mathbf{X}=(X_{1},X_{2})$ has a BDsIW distribution
with the parameter vector $\Psi$ = $\left(  \theta_{1},\theta_{1},\theta
_{1},\zeta\right)  ^{T}$. We will denote this discrete bivariate distribution
by BDsIW$\left(  \theta_{1},\theta_{1},\theta_{1},\zeta\right)  $. If
$\mathbf{X}$ $\sim$ BDsIW$\left(  \theta_{1},\theta_{1},\theta_{1}%
,\zeta\right)  $, then the joint CDF of $\mathbf{X}$ for $x_{1},x_{2}%
\in\mathbf{%
\mathbb{N}
}_{0}$ and for $x_{3}=\min\{x_{1},x_{2}\}$ is given by
\begin{align}
F_{X_{1},X_{2}}(x_{1},x_{2};\Psi)  & =\theta_{1}^{\left(  x_{1}+1\right)
^{-\zeta}}\theta_{2}^{\left(  x_{2}+1\right)  ^{-\zeta}}\theta_{3}^{\left(
x_{3}+1\right)  ^{-\zeta}}\nonumber\\
& =F_{\text{DsIW}}\left(  x_{1};\theta_{1},\zeta\right)  \text{ }%
F_{\text{DsIW}}\left(  x_{2};\theta_{2},\zeta\right)  \text{ }F_{\text{DsIW}%
}\left(  x_{3};\theta_{3},\zeta\right) \nonumber\\
& =\left\{
\begin{array}
[c]{l}%
F_{\text{DsIW}}\left(  x_{1};\theta_{1}\theta_{3},\zeta\right)  \text{
}F_{\text{DsIW}}\left(  x_{2};\theta_{2},\zeta\right)  \ \ \ \ \ \ \ \ \text{;
\ }0<x_{1}<x_{2}<\infty\\
F_{\text{DsIW}}\left(  x_{1};\theta_{1},\zeta\right)  \text{ }F_{\text{DsIW}%
}\left(  x_{2};\theta_{2}\theta_{3},\zeta\right)  \ \ \ \ \ \ \ \ \text{;
\ }0<\text{\ }x_{2}<x_{1}<\infty\\
F_{\text{DsIW}}\left(  x;\theta_{1}\theta_{2}\theta_{3},\zeta\right)  \text{
}\ \ \ \ \ \ \ \ \ \ \ \ \ \ \ \ \ \ \ \ \ \ \ \ \ \text{;\ \ }0<x_{1}%
=x_{2}=x<\infty.
\end{array}
\right. \label{A1}%
\end{align}
The marginal CDFs for BDsIW distribution can be represented as follows%
\begin{equation}
F_{X_{d}}(x_{d};\theta_{d},\theta_{3},\zeta)=P[\max(W_{d},W_{3})\leq
x_{d}]=F_{\text{DsIW}}\left(  x_{d};\theta_{d}\theta_{3},\zeta\right)
.\label{A0}%
\end{equation}
The corresponding joint PMF of $\mathbf{X}$ for $x_{1},x_{2}\in\mathbf{%
\mathbb{N}
}_{0}$ is given by%
\begin{equation}
f_{X_{1},X_{2}}(x_{1},x_{2};\Psi)=\left\{
\begin{array}
[c]{l}%
f_{1}(x_{1},x_{2};\Psi)\text{ }\ \ \ \ \ \ \ \ \ \ \ \ \ \text{; \ }%
0<x_{1}<x_{2}<\infty\\
f_{2}(x_{1},x_{2};\Psi)\text{ }\ \ \ \ \ \ \ \ \ \ \ \ \ \text{;
\ }0<\text{\ }x_{2}<x_{1}<\infty\\
f_{3}(x;\Psi)\ \ \ \ \ \ \ \ \ \ \ \ \ \ \ \ \ \ \ \ \ \text{;\ \ }%
0<x_{1}=x_{2}=x<\infty,
\end{array}
\right. \label{A2}%
\end{equation}
where%
\begin{align*}
f_{1}(x_{1},x_{2};\Psi)  & =f_{\text{DsIW}}\left(  x_{1};\theta_{1}\theta
_{3},\zeta\right)  f_{\text{DsIW}}\left(  x_{2};\theta_{2},\zeta\right)  ,\\
f_{2}(x_{1},x_{2};\Psi)  & =f_{\text{DsIW}}\left(  x_{1};\theta_{1}%
,\zeta\right)  \text{ }f_{\text{DsIW}}\left(  x_{2};\theta_{2}\theta_{3}%
,\zeta\right)  ,\\
f_{3}(x;\Psi)  & =F_{\text{DsIW}}\left(  x;\theta_{2},\zeta\right)
f_{\text{DsIW}}\left(  x;\theta_{1}\theta_{3},\zeta\right)  -F_{\text{DsIW}%
}\left(  x-1;\theta_{2}\theta_{3},\zeta\right)  f_{\text{DsIW}}\left(
x;\theta_{1},\zeta\right)  .
\end{align*}
The expressions $f_{1}(x_{1},x_{2};\Psi)$, $f_{2}(x_{1},x_{2};\Psi)$ and
$f_{3}(x;\Psi)$ for $x_{1},x_{2}\in\mathbf{%
\mathbb{N}
}_{0}$ can be easily obtained by using the relation%
\begin{equation}
f_{X_{1},X_{2}}(x_{1},x_{2};\Psi)=F(x_{1},x_{2};\Psi)-F(x_{1}-1,x_{2}%
;\Psi)-F(x_{1},x_{2}-1;\Psi)+F(x_{1}-1,x_{2}-1;\Psi).
\end{equation}
Figure 1 shows the plots of the joint PMF of BDsIW distribution for different
parameter values.

\begin{center}%
\[%
{\includegraphics[
height=1.9579in,
width=2.0349in
]%
{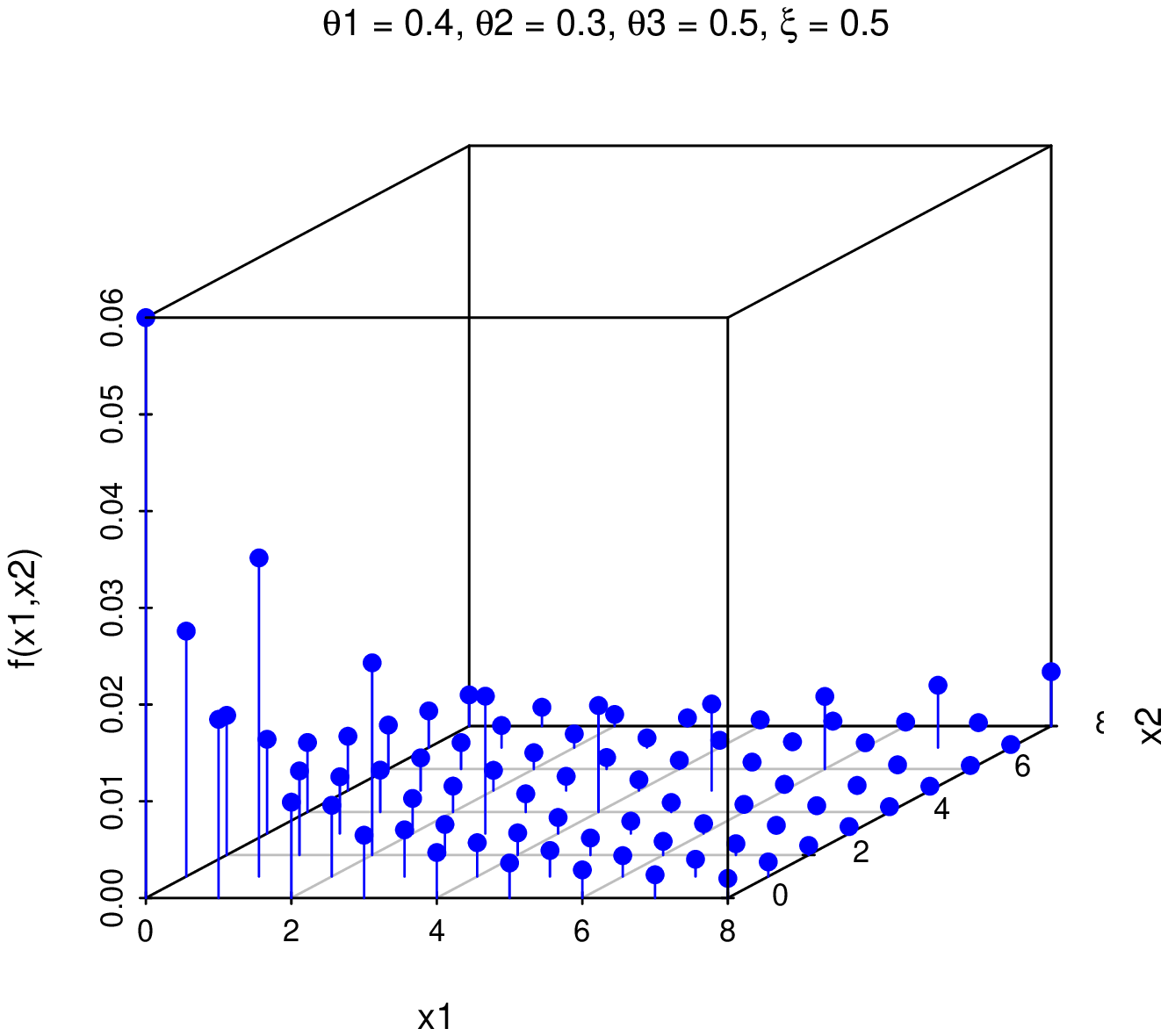}%
}
{\includegraphics[
height=1.9579in,
width=2.0349in
]%
{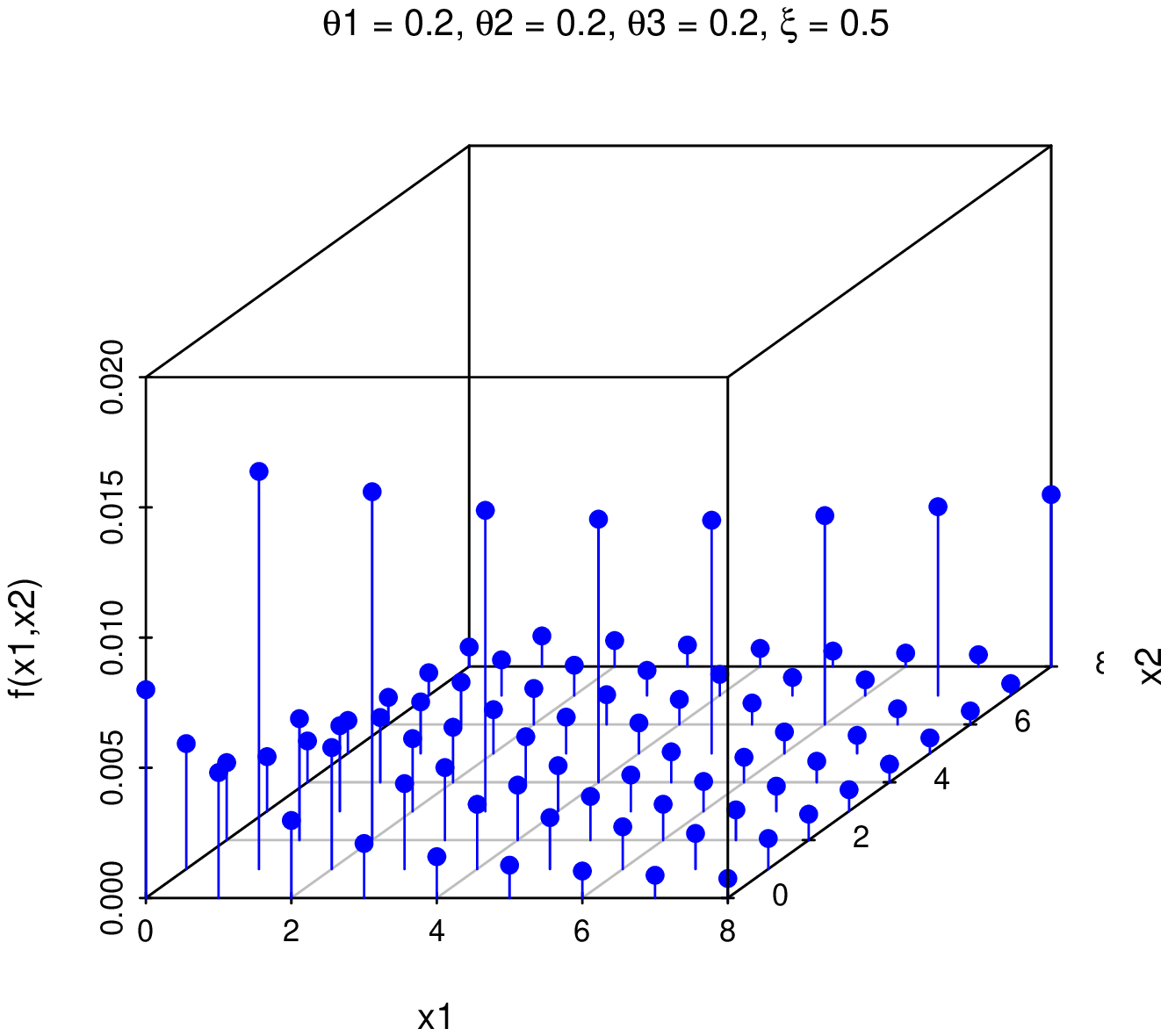}%
}
{\includegraphics[
height=1.9579in,
width=2.0349in
]%
{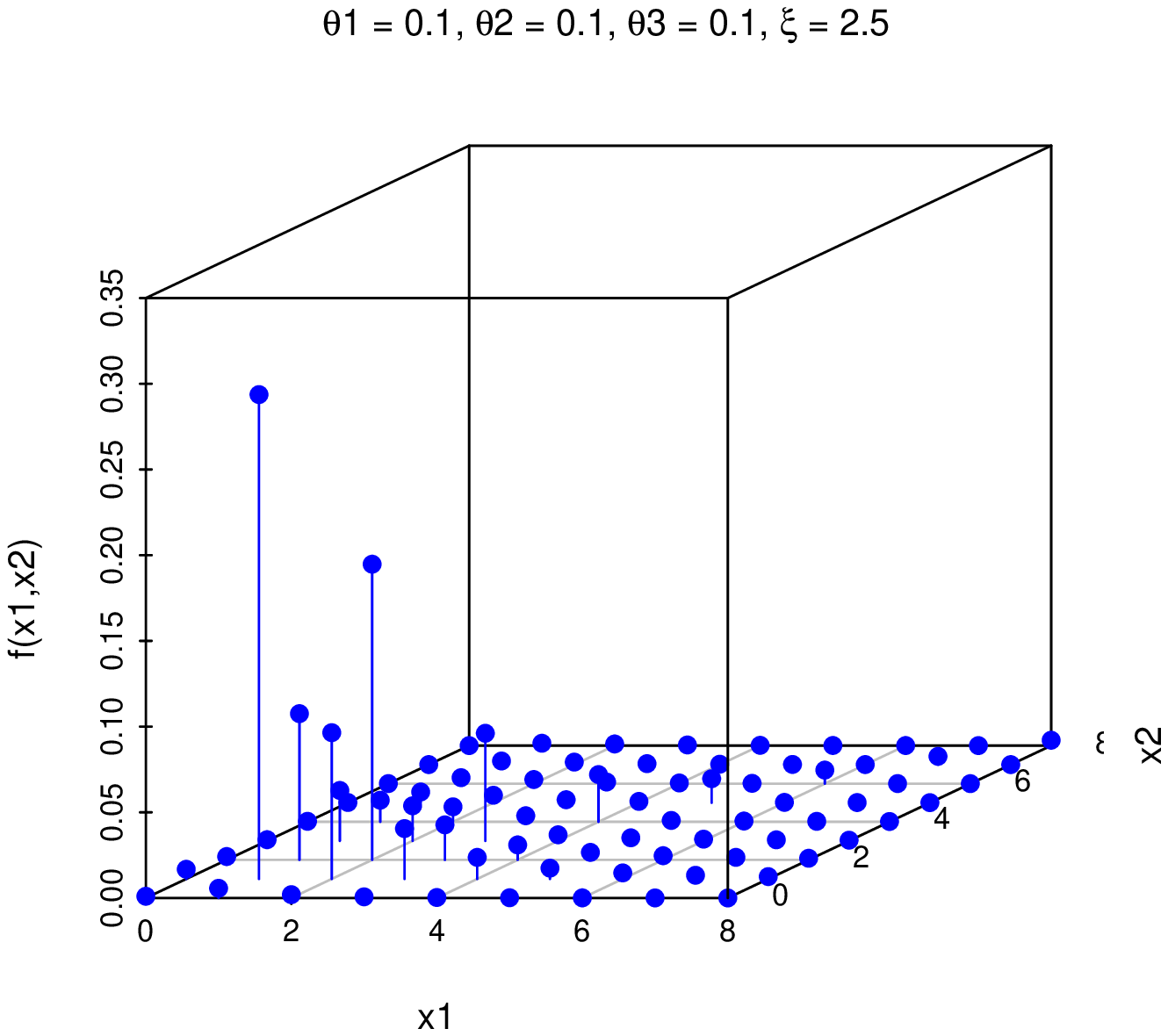}%
}
\]
Figure 1. The scatter plots of the joint PMF for different parameter values.
\end{center}

From Figure 1, it is clear that the joint PMF can take different shapes
depending on the model parameter values. Assume $\mathbf{X}$ $\sim$
BDsIW$\left(  \theta_{1},\theta_{2},\theta_{3},\zeta\right)  $,$\ $then the
joint reliability function of $\mathbf{X}\ $can be expressed as%
\begin{align}
R_{X_{1},X_{2}}(x_{1},x_{2};\Psi)  & =1-F_{X_{1}}(x_{1};\theta_{1}\theta
_{3},\zeta)\ -F_{X_{2}}(x_{2};\theta_{2}\theta_{3},\zeta)+F_{X_{1},X_{2}%
}(x_{1},x_{2};\Psi)\nonumber\\
& =\left\{
\begin{array}
[c]{l}%
R_{1}(x_{1},x_{2};\Psi)\text{ }\ \ \ \ \ \ \ \ \ \ \ \ \ \text{; \ }%
0<x_{1}<x_{2}<\infty\\
R_{2}(x_{1},x_{2};\Psi)\text{ }\ \ \ \ \ \ \ \ \ \ \ \ \ \text{;
\ }0<\text{\ }x_{2}<x_{1}<\infty\\
R_{3}(x;\Psi)\ \ \ \ \ \ \ \ \ \ \ \ \ \ \ \ \ \ \ \ \ \text{;\ \ }%
0<x_{1}=x_{2}=x<\infty,
\end{array}
\right. \label{A3}%
\end{align}
where%
\begin{align*}
R_{1}(x_{1},x_{2};\Psi)\text{ }  & =1-F_{\text{DsIW}}(x_{1};\theta_{1}%
\theta_{3},\zeta)-F_{\text{DsIW}}(x_{2};\theta_{2}\theta_{3},\zeta
)+F_{\text{DsIW}}\left(  x_{1};\theta_{1}\theta_{3},\zeta\right)  \text{
}F_{\text{DsIW}}\left(  x_{2};\theta_{2},\zeta\right)  ,\\
R_{2}(x_{1},x_{2};\Psi)\text{ }  & =1-F_{\text{DsIW}}(x_{1};\theta_{1}%
\theta_{3},\zeta)-F_{\text{DsIW}}(x_{2};\theta_{2}\theta_{3},\zeta
)+F_{\text{DsIW}}\left(  x_{1};\theta_{1},\zeta\right)  \text{ }%
F_{\text{DsIW}}\left(  x_{2};\theta_{2}\theta_{3},\zeta\right)  ,
\end{align*}
and
\[
R_{3}(x;\Psi)\ \text{ }=1-F_{\text{DsIW}}(x;\theta_{1}\theta_{3}%
,\zeta)-F_{\text{DsIW}}(x;\theta_{2}\theta_{3},\zeta)+F_{\text{DsIW}}\left(
x;\theta_{1}\theta_{2}\theta_{3},\zeta\right)  .
\]
Moreover, the bivariate hazard rate function (BHRF) of $\mathbf{X}$ can be
represented as%
\begin{equation}
r_{X_{1},X_{2}}(x_{1},x_{2};\Psi)=\left\{
\begin{array}
[c]{l}%
r_{1}(x_{1},x_{2};\Psi)\text{ }\ \ \ \ \ \ \ \ \ \ \ \ \ \text{; \ }%
0<x_{1}<x_{2}<\infty\\
r_{2}(x_{1},x_{2};\Psi)\text{ }\ \ \ \ \ \ \ \ \ \ \ \ \ \text{;
\ }0<\text{\ }x_{2}<x_{1}<\infty\\
r_{3}(x;\Psi)\ \ \ \ \ \ \ \ \ \ \ \ \ \ \ \ \ \ \ \ \ \text{;\ \ }%
0<x_{1}=x_{2}=x<\infty,
\end{array}
\right.
\end{equation}
where $r_{j}(.;\Psi)$ =$\frac{f_{j}(.;\Psi)}{R_{j}(.;\Psi)};j=1,2,3.$ Figure 2
shows the plots of the BHRF of BDsIW distribution for different parameter values.

\begin{center}%
\[%
{\includegraphics[
height=1.9579in,
width=2.0349in
]%
{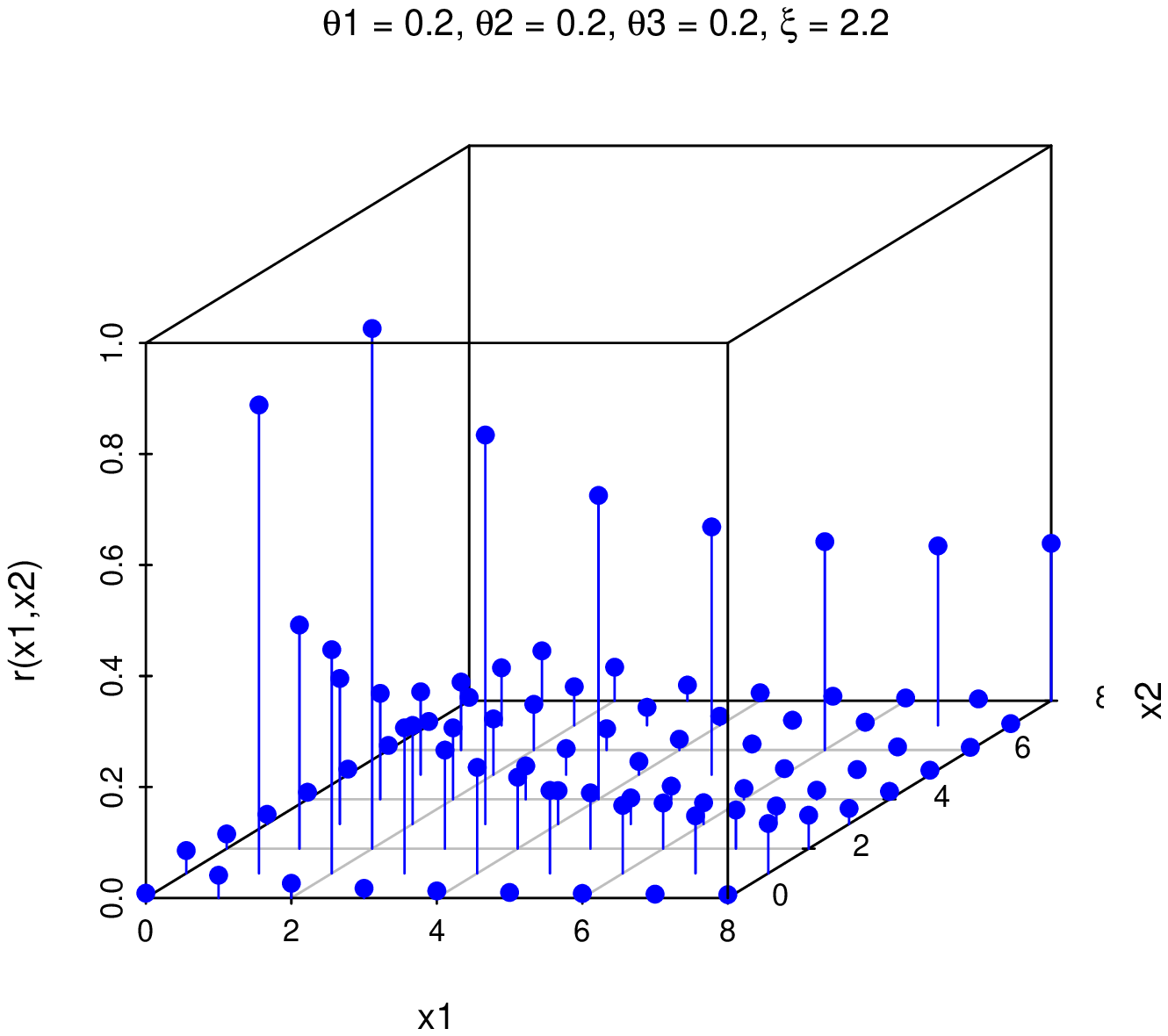}%
}
{\includegraphics[
height=1.9579in,
width=2.0349in
]%
{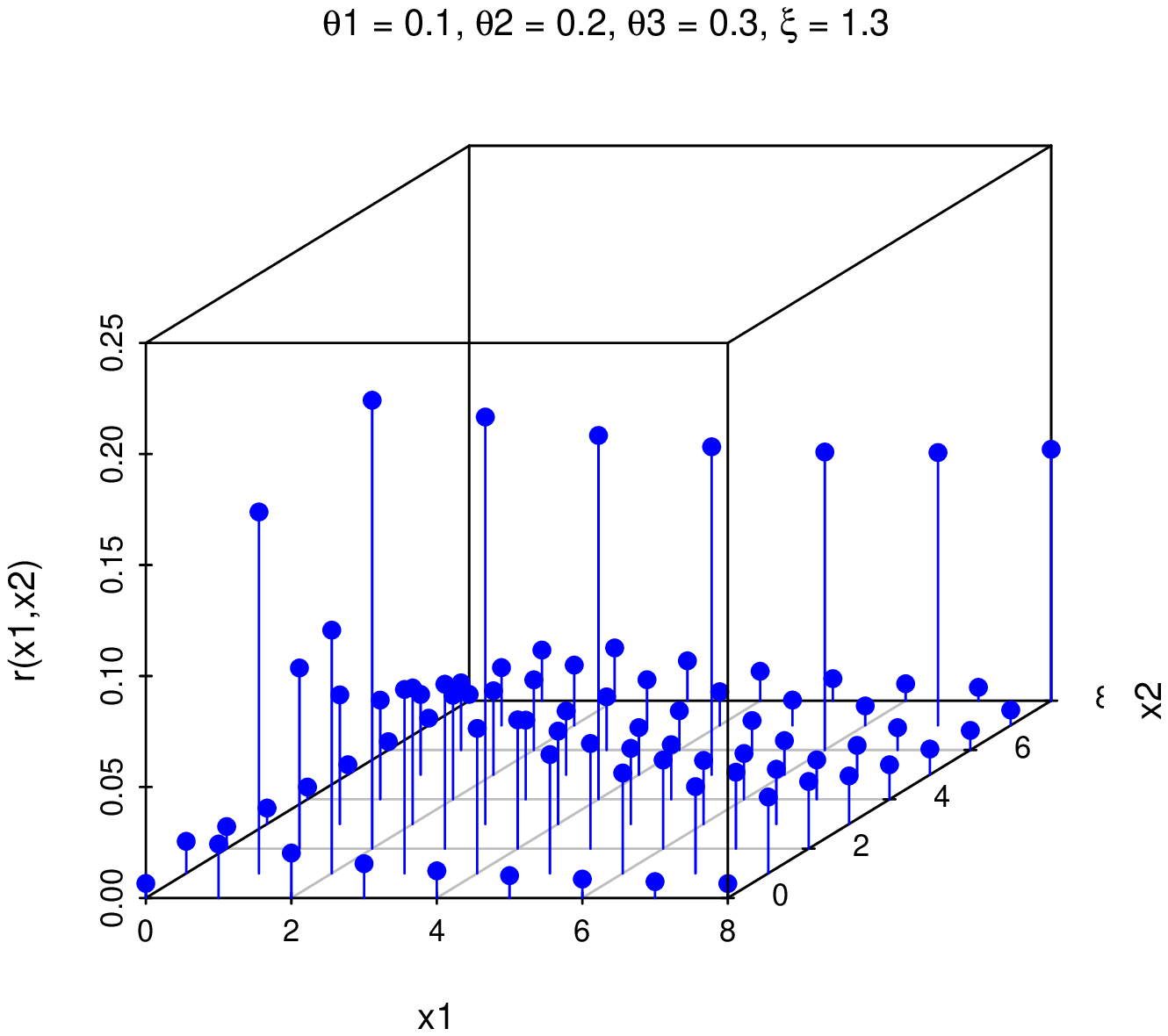}%
}
{\includegraphics[
height=1.9579in,
width=2.0349in
]%
{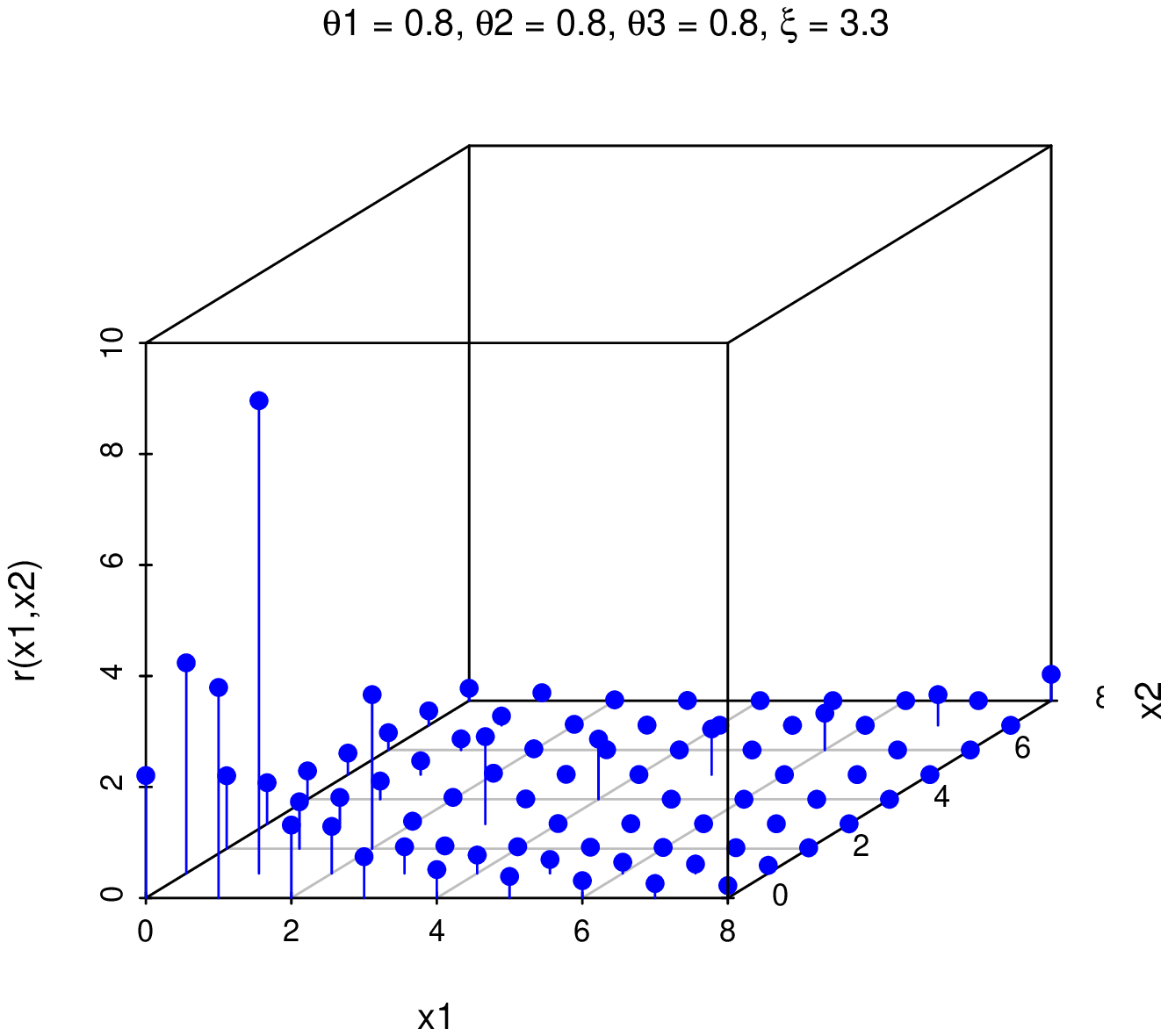}%
}
\]
Figure 2. The scatter plots of the BHRF for different parameter values.
\end{center}

From Figure 2, it is clear that the joint BHRF can take different shapes
depending on the parameter values. Assume $X_{1}$ $<$ $X_{2}$, then the HRF of
the conditional distribution $X_{1}$ given $X_{2}>x_{2}$ is given by%
\begin{equation}
r^{\ast}(X_{1}|X_{2}>x_{2})=\frac{\zeta\left(  x_{1}+1\right)  ^{-\zeta-1}%
}{R_{1}(x_{1},x_{2};\Psi)}\left\{  F_{\text{DsIW}}\left(  x_{2};\theta
_{2},\zeta\right)  -1\right\}  F_{\text{DsIW}}\left(  x_{1};\theta_{1}%
\theta_{3},\zeta\right)  \ln\left(  \theta_{1}\theta_{3}\right)  ,
\end{equation}
and the HRF of the conditional distribution $X_{2}$ given $X_{1}>x_{1}$ is
given by%
\begin{equation}
r^{\ast}(X_{2}|X_{1}>x_{1})=\frac{\zeta\left(  x_{2}+1\right)  ^{-\zeta
-1}F_{\text{DsIW}}\left(  x_{2};\theta_{2},\zeta\right)  }{R_{1}(x_{1}%
,x_{2};\Psi)}\left\{  F_{\text{DsIW}}\left(  x_{1};\theta_{1}\theta_{3}%
,\zeta\right)  \ln\left(  \theta_{2}\right)  -F_{\text{DsIW}}\left(
x_{2};\theta_{3},\zeta\right)  \ln\left(  \theta_{2}\theta_{3}\right)
\right\}  .
\end{equation}
Similarly, when $X_{2}$ $<$ $X_{1}$, then%
\begin{equation}
r^{\ast\ast}(X_{1}|X_{2}>x_{2})=\frac{\zeta\left(  x_{1}+1\right)  ^{-\zeta
-1}F_{\text{DsIW}}\left(  x_{1};\theta_{1},\zeta\right)  }{R_{2}(x_{1}%
,x_{2};\Psi)}\left\{  F_{\text{DsIW}}\left(  x_{2};\theta_{2}\theta_{3}%
,\zeta\right)  \ln\left(  \theta_{1}\right)  -F_{\text{DsIW}}\left(
x_{1};\theta_{3},\zeta\right)  \ln\left(  \theta_{1}\theta_{3}\right)
\right\}  ,
\end{equation}
and%
\begin{equation}
r^{\ast\ast}(X_{2}|X_{1}>x_{1})=\frac{\zeta\left(  x_{2}+1\right)  ^{-\zeta
-1}}{R_{2}(x_{1},x_{2};\Psi)}\left\{  F_{\text{DsIW}}\left(  x_{1};\theta
_{1},\zeta\right)  -1\right\}  F_{\text{DsIW}}\left(  x_{2};\theta_{2}%
\theta_{3},\zeta\right)  \ln\left(  \theta_{2}\theta_{3}\right)  .
\end{equation}

On the other hand, assume a parallel system contains 2-component. Then, we can
defined the BHRF as a vector which is useful to measure the total life span of
a 2-component as follows%
\begin{equation}
r(\underline{x})=(r(x),r_{12}(x_{1}|x_{2}),r_{21}(x_{2}|x_{1})),
\end{equation}
where $r(x)$ gives the HRF of the system using the information that the
2-component has survived beyond $x$, $r_{12}(x_{1}|x_{2})$ gives the HRF span
of the first component given that it has survived to an age $x_{1}$ and the
other has failed at $x_{1}$. Similar argument holds for $r_{21}(x_{2}|x_{1})$,
(see Cox (1972)). If $\mathbf{X}$ $\sim$ BDsIW$\left(  \theta_{1},\theta
_{1},\theta_{1},\zeta\right)  $, then%
\begin{align*}
r(x)|_{X=\min\left(  X_{1},X_{2}\right)  }  & =\frac{F_{\text{DsIW}}\left(
x-1;\theta_{3},\zeta\right)  }{R_{3}(x;\Psi)}\left[  -F_{\text{DsIW}}\left(
x-1;\theta_{1},\zeta\right)  -F_{\text{DsIW}}\left(  x-1;\theta_{2}%
,\zeta\right)  +F_{\text{DsIW}}\left(  x-1;\theta_{1}\theta_{2},\zeta\right)
\right] \\
& +\frac{F_{\text{DsIW}}\left(  x;\theta_{3},\zeta\right)  }{R_{3}(x;\Psi
)}\left[  F_{\text{DsIW}}\left(  x;\theta_{1},\zeta\right)  +F_{\text{DsIW}%
}\left(  x;\theta_{2},\zeta\right)  -F_{\text{DsIW}}\left(  x;\theta_{1}%
\theta_{2},\zeta\right)  \right]  ,
\end{align*}%
\[
r_{12}(x_{1}|x_{2})|_{X_{1}>X_{2}}=\zeta\left(  x_{1}+1\right)  ^{-\zeta-1}
\left[  1-F_{\text{DsIW}}\left(  x_{1};\theta_{1},\zeta\right)  \right]
^{-1}\ln\left(  \theta_{1}\right)  ,
\]
and%
\[
r_{21}(x_{2}|x_{1})|_{X_{1}<X_{2}}=\zeta\left(  x_{2}+1\right)  ^{-\zeta-1}
\left[  1-F_{\text{DsIW}}\left(  x_{2};\theta_{2},\zeta\right)  \right]
^{-1}\ln\left(  \theta_{2}\right)  .
\]

The following shock model and maintenance model interpretations can be
provided for BDsIW distribution.

\begin{enumerate}
\item Shock model: Consider a system has 2-component, and it is assumed that
the amount of shocks is measured in a discrete unit. Each component is
subjected to individual shocks say $W_{1}$ and $W_{2}$ respectively. The
system faces an overall shock $W_{3}$, which is transmitted to both the
component equally, independent of their individual shocks. So, the observed
shocks at the 2- component are $X_{1}=\max(W_{1},W_{3})$ and $X_{2}=\max
(W_{2},W_{3})$ respectively.

\item Maintenance model: Consider a system has 2-component, and it is assumed
that each component has been maintained independently and also there is an
overall maintenance. Due to component maintenance, assume the lifetime of the
individual component is increased by $W_{i}$ amount and because of the overall
maintenance, the lifetime of each component is increased by $W_{3}$ amount.
Here, $W_{1}$, $W_{2}$ and $W_{3}$ are all measured in a discrete unit. So,
the increased lifetimes of the 2-component are $X_{1}=\max(W_{1},W_{3})$ and
$X_{2}=\max(W_{2},W_{3})$ respectively.
\end{enumerate}

\subsection{Some statistical properties}

Assume $\mathbf{X}$ $\sim$ BDsIW$\left(  \theta_{1},\theta_{2},\theta
_{3},\zeta\right)  $, then $X_{1}$and$\ X_{2}$ are positive quadrant dependent
(PQD) where%
\begin{equation}
F_{X_{1},X_{2}}(x_{1},x_{2};\Psi)\geq F_{X_{1}}(x_{1};\theta_{1},\theta
_{3},\zeta)F_{X_{2}}(x_{2};\theta_{2},\theta_{3},\zeta).
\end{equation}
Furthermore, for every pair of increasing functions $f_{X_{1}}(.)$ and
$f_{X_{2}}(.)$, we get $Cov\left\{  f_{X_{1}}(X_{1}),f_{X_{2}}(X_{2})\right\}
\geq0$; see for example Nelsen (2006).\ Let us recall that, the function
$k(u,v):R\times R\rightarrow R$, is said to have the total positivity of order
two ($TP-O_{2}$) property if $k(u,v)$ satisfies%
\begin{equation}
k(u_{1},v_{1})k(u_{2},v_{2)}\geq k(u_{2},v_{1})k(u_{1},v_{2}),
\end{equation}
for all $u_{1},v_{1},u_{2},v_{2}\in R$. It is consider a very strong and an
important property in lifetime testing, see for example Hu et al. (2003).
Assume $x_{11},x_{21},x_{12},x_{22}\in\mathbf{%
\mathbb{N}
}_{0}$ and $x_{11}<x_{21}<x_{12}<x_{22}$ from $\mathbf{X}$ $\sim$
BDsIW$\left(  \theta_{1},\theta_{2},\theta_{3},\zeta\right)  $, then the joint
reliability function of $\mathbf{X}$ satisfies the $TP-O_{2}$ property where
\begin{equation}
\frac{R_{X_{1},X_{2}}(x_{11},x_{21})R_{X_{1},X_{2}}(x_{12},x_{22})}%
{R_{X_{1},X_{2}}(x_{12},x_{21})R_{X_{1},X_{2}}(x_{11},x_{22})}\geq1.
\end{equation}
Similarly, when $x_{11}=x_{21}<x_{12}<x_{22},$ $x_{21}<x_{11}<x_{12}<x_{22}$
etc. Now, we present some statistical properties of the proposed model in a
form of results.

\textbf{Result 1.} If the bivariate vector $\mathbf{X}=(X_{1},X_{2})$ has the
BDsIW$\left(  \theta_{1},\theta_{2},\theta_{3},\zeta\right)  $, then

\begin{enumerate}
\item $\max\{X_{1},X_{2}\}\sim$ DsIW$\left(  \theta_{1}\theta_{2}\theta
_{3},\zeta\right)  .$

\item The stress-strenght probability is given by%
\begin{equation}
P[X_{1}<X_{2}]=\sum_{x=0}^{\infty}\left(  \theta_{1}\theta_{3}\right)
^{\left(  x+1\right)  ^{-\zeta}}\left[  \left(  \theta_{2}\theta_{3}\right)
^{\left(  x+1\right)  ^{-\zeta}}-\left(  \theta_{2}\theta_{3}\right)
^{x^{-\zeta}}\right]  .
\end{equation}

\item The median of $X_{1}$ and $X_{2}$ is given by%
\begin{equation}
M_{X_{d}}=\left\{  \log\frac{\theta_{d}\theta_{3}}{U}\right\}  ^{\frac
{1}{\zeta}}-1;\ d=1,2,
\end{equation}
where $U$ \ has a uniform $U(0,1)$ distribution.

\item The coefficient of median correlation between $X_{1}$ and $X_{2}$ is
given by%
\begin{equation}
M_{X_{1},X_{2}}=\left\{
\begin{array}
[c]{l}%
4F_{\text{DsIW}}\left(  M_{X_{1}};\theta_{1}\theta_{3},\zeta\right)  \text{
}F_{\text{DsIW}}\left(  M_{X_{2}};\theta_{2},\zeta\right)
-1\ \ \ \ \ \ \ \ \text{; \ }x_{1}<x_{2}\\
4F_{\text{DsIW}}\left(  M_{X_{1}};\theta_{1},\zeta\right)  \text{
}F_{\text{DsIW}}\left(  M_{X_{2}};\theta_{2}\theta_{3},\zeta\right)
-1\ \ \ \ \ \ \ \ \text{; \ }x_{2}\leq x_{1}.
\end{array}
\right.
\end{equation}

\item The conditional PMF of $X_{1}$ given $X_{2}=x_{2}$, is given by%
\begin{equation}
f_{X_{1}\mid X_{2}=x_{2}}(x_{1}\mid x_{2})=\left\{
\begin{array}
[c]{l}%
f_{X_{1}\mid X_{2}=x_{2}}^{(1)}(x_{1}\mid x_{2})\ \ \ \ \text{if }%
\ 0<x_{1}<x_{2}<\infty\\
f_{X_{1}\mid X_{2}=x_{2}}^{(2)}(x_{1}\mid x_{2})\ \ \ \ \text{if \ }%
0<x_{2}<x_{1}<\infty\\
f_{X_{1}\mid X_{2}=x_{2}}^{(3)}(x_{1}\mid x)\ \ \ \ \text{if\ \ }0<x_{1}%
=x_{2}=x<\infty,
\end{array}
\right.
\end{equation}
where%
\begin{align*}
f_{X_{1}\mid X_{2}=x_{2}}^{(1)}(x_{1}  & \mid x_{2})=\frac{f_{\text{DsIW}%
}\left(  x_{1};\theta_{1}\theta_{3},\zeta\right)  f_{\text{DsIW}}\left(
x_{2};\theta_{2},\zeta\right)  }{f_{\text{DsIW}}\left(  x_{2};\theta_{2}%
\theta_{3},\zeta\right)  },\\
f_{X_{1}\mid X_{2}=x_{2}}^{(2)}(x_{1}  & \mid x_{2})=f_{\text{DsIW}}\left(
x_{1};\theta_{1},\zeta\right)  ,
\end{align*}
and%
\[
f_{X_{1}\mid X_{2}=x_{2}}^{(3)}(x_{1}\mid x)\ =\frac{F_{\text{DsIW}}\left(
x;\theta_{2},\zeta\right)  f_{\text{DsIW}}\left(  x;\theta_{1}\theta_{3}%
,\zeta\right)  -F_{\text{DsIW}}\left(  x-1;\theta_{2}\theta_{3},\zeta\right)
f_{\text{DsIW}}\left(  x;\theta_{1},\zeta\right)  }{f_{\text{DsIW}}\left(
x;\theta_{2}\theta_{3},\zeta\right)  }.
\]

\item The conditional CDF of $X_{1}$ given $X_{2}\leq x_{2}$, is given by%
\begin{equation}
F_{X_{1}\mid X_{2}=x_{2}}(x_{1}\mid x_{2})=\left\{
\begin{array}
[c]{l}%
\frac{F_{\text{DsIW}}\left(  x_{1};\theta_{1}\theta_{3},\zeta\right)  \text{
}}{F_{\text{DsIW}}\left(  x_{2};\theta_{3},\zeta\right)  }\ \ \ \ \ \ \text{if
}\ 0<x_{1}<x_{2}<\infty\\
F_{\text{DsIW}}\left(  x_{1};\theta_{1},\zeta\right)  \ \ \ \ \text{if
\ }0<x_{2}<x_{1}<\infty\\
F_{\text{DsIW}}\left(  x;\theta_{1},\zeta\right)  \ \ \ \ \ \ \text{if\ \ }%
0<x_{1}=x_{2}=x<\infty.
\end{array}
\right.
\end{equation}

\end{enumerate}

\textbf{Proof}. The proofs are quite standard and the details are avoided.$\ $

\textbf{Result 2. }Assume $(X_{i1},X_{i2})$ $\sim$ BDsIW$\left(  \theta
_{i1},\theta_{i2},\theta_{i3},\zeta\right)  $ for $i=1,2,...,n$ and they are
independently distributed. If%
\[
Z_{s}=\max\left(  X_{1s},X_{2s},...,X_{ns}\right)  ;s=1,2\Longrightarrow
(X_{i1},X_{i2})\sim\text{BDsIW}\left(
{\textstyle\prod\limits_{i=1}^{n}}
\theta_{i1},%
{\textstyle\prod\limits_{i=1}^{n}}
\theta_{i2},%
{\textstyle\prod\limits_{i=1}^{n}}
\theta_{i3},\zeta\right)  .
\]

\textbf{Proof. }It is easy to proof that using the joint CDF.

\textbf{Result 3.} If the bivariate vector $\mathbf{X}$ $\sim$ BDsIW$\left(
\theta_{1},\theta_{2},\theta_{3},\zeta\right)  $, then the joint probability
generating function (PGF) of $X_{1}$ and$\ X_{2}$ can be written as infinite
mixtures,%
\begin{align}
G_{X_{1},X_{2}}\left(  y_{1},y_{2}\right)   & =\sum_{j=0}^{\infty}\sum
_{i=0}^{j-1}\left[  \left(  \theta_{1}\theta_{3}\right)  ^{\left(  i+1\right)
^{-\zeta}}-\left(  \theta_{1}\theta_{3}\right)  ^{i^{-\zeta}}\right]  \left[
\theta_{2}^{\left(  j+1\right)  ^{-\zeta}}-\theta_{2}^{j^{-\zeta}}\right]
y_{1}^{i}y_{2}^{j}\nonumber\\
& +\sum_{j=0}^{\infty}\sum_{i=j+1}^{\infty}\left[  \theta_{1}^{\left(
i+1\right)  ^{-\zeta}}-\theta_{1}^{i^{-\zeta}}\right]  \left[  \left(
\theta_{2}\theta_{3}\right)  ^{\left(  j+1\right)  ^{-\zeta}}-\left(
\theta_{2}\theta_{3}\right)  ^{j^{-\zeta}}\right]  y_{1}^{i}y_{2}%
^{j}\nonumber\\
& +\sum_{i=0}^{\infty}\theta_{2}^{\left(  i+1\right)  ^{-\zeta}}\left[
\left(  \theta_{1}\theta_{3}\right)  ^{\left(  i+1\right)  ^{-\zeta}}-\left(
\theta_{1}\theta_{3}\right)  ^{i^{-\zeta}}\right]  y_{1}^{i}y_{2}%
^{i}\nonumber\\
& -\sum_{i=0}^{\infty}\left(  \theta_{2}\theta_{3}\right)  ^{\left(
i+1\right)  ^{-\zeta}}\left[  \theta_{1}^{\left(  i+1\right)  ^{-\zeta}%
}-\theta_{1}^{i^{-\zeta}}\right]  y_{1}^{i}y_{2}^{i};\ \ \ \ \left\vert
y_{1}\right\vert ,\left\vert y_{2}\right\vert <1.
\end{align}

\textbf{Proof}. The proof can be easily obtained by using the fact that%
\[
G_{X_{1},X_{2}}\left(  y_{1},y_{2}\right)  =E\left(  y_{1}^{X_{1}}y_{2}%
^{X_{2}}\right)  =\sum_{j,i=0}^{\infty}P\left[  X_{1}=i,X_{2}=j\right]
y_{1}^{i}y_{2}^{j}.
\]
Hence, different moments and product moments of BDsIW distribution can be
obtained, as infinite series, using the joint PGF.

\section{Statistical Inference}

\subsection{Maximum likelihood estimation (MLE)}

In this section, we use the maximum likelihood method to estimate the unknown
parameters $\theta_{1},\theta_{2},\theta_{3}$ and $\zeta$ of\ BDsIW
distribution. Suppose that, we have a sample of size $n$, of the form
$\left\{  (x_{11},x_{21}),(x_{12},x_{22}),...,(x_{1n},x_{2n})\right\}  $ from
BDsIW distribution. We use the following notations: $I_{1}=\{x_{1j}<x_{2j}\},
$ $I_{2}=\{x_{2j}<x_{1j}\},$ $I_{3}=\{x_{1j}=x_{2j}=x_{j}\},$ $I=I_{1}\cup
I_{2}\cup I_{3},$ $\left\vert I_{1}\right\vert =n_{1},$ $\left\vert
I_{2}\right\vert =n_{2},$ $\left\vert I_{3}\right\vert =n_{3}$ and
$n=\sum_{k=1}^{3}n_{k}.$ Based on the observations, the likelihood function is
given by%

\begin{equation}
l(\Psi)=\underset{j=1}{\overset{n_{1}}{\prod}}f_{1}(x_{1j},x_{2j}%
)\underset{j=1}{\overset{n_{2}}{\prod}}f_{2}(x_{1j},x_{2j})\underset
{j=1}{\overset{n_{3}}{\prod}}f_{3}(x_{j}).
\end{equation}
The log-likelihood function becomes%

\begin{align}
L(\Psi)  & =\overset{n_{1}}{\underset{j=1}{\sum}}\ln\left(  \Phi_{1}%
(x_{1j};\theta_{1}\theta_{3},\zeta)\right)  +\overset{n_{1}}{\underset
{j=1}{\sum}}\ln\left(  \Phi_{1}(x_{2j};\theta_{2},\zeta)\right)
+\overset{n_{2}}{\underset{j=1}{\sum}}\ln\left(  \Phi_{1}(x_{1j};\theta
_{1},\zeta)\right)  +\overset{n_{2}}{\underset{j=1}{\sum}}\ln\left(  \Phi
_{1}(x_{2j};\theta_{2}\theta_{3},\zeta)\right) \nonumber\\
& +\overset{n_{3}}{\underset{j=1}{\sum}}\ln\left(  \left[  \theta_{2}\right]
^{(x_{j}+1)^{-\zeta}}\Phi_{1}(x_{j};\theta_{1}\theta_{3},\zeta)-\left[
\theta_{2}\theta_{3}\right]  ^{(x_{j})^{-\zeta}}\Phi_{1}(x_{j};\theta
_{1},\zeta)\right)  ,\label{1.36}%
\end{align}
where $\Phi_{1}(x;\theta,\zeta)=\theta^{(x+1)^{\zeta}}-\theta^{x^{\zeta}}.$The
MLEs of the model parameters can be obtained by computing the first partial
derivatives of Equation (\ref{1.36}) with respect to $\theta_{1},\theta
_{2},\theta_{3}$ and $\zeta$ and then putting the results equal zeros. We get
the likelihood equations as in the following form%
\begin{align}
\frac{\partial L}{\partial\theta_{1}}  & =\overset{n_{1}}{\underset{j=1}{\sum
}}\frac{\theta_{3}\Phi_{2}(x_{1j}+1;\theta_{1}\theta_{3},\zeta)-\theta_{3}%
\Phi_{2}(x_{1j};\theta_{1}\theta_{3},\zeta)}{\Phi_{1}(x_{1j};\theta_{1}%
\theta_{3},\zeta)}+\overset{n_{2}}{\underset{j=1}{\sum}}\frac{\Phi_{2}%
(x_{1j}+1;\theta,\zeta)-\Phi_{2}(x_{1j};\theta,\zeta)}{\Phi_{1}(x_{1j}%
;\theta_{1},\zeta)}+\nonumber\\
& \overset{n_{3}}{\underset{j=1}{\sum}}\frac{\theta_{3}\theta_{2}%
^{(x_{j}+1)^{-\zeta}}\left[  \Phi_{2}(x_{j}+1;\theta_{1}\theta_{3},\zeta
)-\Phi_{2}(x_{j};\theta_{1}\theta_{3},\zeta)\right]  -\left(  \theta_{2}%
\theta_{3}\right)  ^{(x_{j})^{-\zeta}}\left[  \Phi_{2}(x_{j}+1;\theta
_{1},\zeta)-\Phi_{2}(x_{j};\theta_{1},\zeta)\right]  }{\theta_{2}%
^{(x_{j}+1)^{-\zeta}}\Phi_{1}(x_{j};\theta_{1}\theta_{3},\zeta)-\left(
\theta_{2}\theta_{3}\right)  ^{(x_{j})^{-\zeta}}\Phi_{1}(x_{j};\theta
_{1},\zeta)},\label{A5}%
\end{align}%
\begin{align}
\frac{\partial L}{\partial\theta_{2}}  & =\overset{n_{1}}{\underset{j=1}{\sum
}}\frac{\Phi_{2}(x_{2j}+1;\theta_{2},\zeta)-\Phi_{2}(x_{2j};\theta_{2},\zeta
)}{\Phi_{1}(x_{2j};\theta_{2},\zeta)}+\overset{n_{2}}{\underset{j=1}{\sum}%
}\frac{\theta_{3}\Phi_{2}(x_{2j}+1;\theta_{2}\theta_{3},\zeta)-\theta_{3}%
\Phi_{2}(x_{2j};\theta_{2}\theta_{3},\zeta)}{\Phi_{1}(x_{2j};\theta_{1}%
,\zeta)}\nonumber\\
& \overset{n_{3}}{\underset{j=1}{+\sum}}\frac{\Phi_{2}(x_{j}+1;\theta
_{2},\zeta)\Phi_{1}(x_{j};\theta_{1}\theta_{3},\zeta)-\theta_{3}\Phi_{2}%
(x_{j};\theta_{2}\theta_{3},\zeta)\Phi_{2}(x_{j};\theta_{1},\zeta)}{\theta
_{2}^{(x_{j}+1)^{-\zeta}}\Phi_{1}(x_{j};\theta_{1}\theta_{3},\zeta)-\left(
\theta_{2}\theta_{3}\right)  ^{(x_{j})^{-\zeta}}\Phi_{1}(x_{j};\theta
_{1},\zeta)},\label{A6}%
\end{align}%
\begin{align}
\frac{\partial L}{\partial\theta_{3}}  & =\overset{n_{1}}{\underset{j=1}{\sum
}}\frac{\theta_{1}\left[  \Phi_{2}(x_{1j}+1;\theta_{1}\theta_{3},\zeta
)-\Phi_{2}(x_{1j};\theta_{1}\theta_{3},\zeta)\right]  }{\Phi_{1}(x_{1j}%
;\theta_{1}\theta_{3},\zeta)}+\overset{n_{2}}{\underset{j=1}{\sum}}%
\frac{\theta_{2}\left[  \Phi_{2}(x_{2j}+1;\theta_{2}\theta_{3},\zeta)-\Phi
_{2}(x_{2j};\theta_{2}\theta_{3},\zeta)\right]  }{\Phi_{1}(x_{2j};\theta
_{2}\theta_{3},\zeta)}+\nonumber\\
& \overset{n_{3}}{\underset{j=1}{\sum}}\frac{\theta_{1}\theta_{2}%
^{(x_{j}+1)^{-\zeta}}\left[  \Phi_{2}(x_{j}+1;\theta_{1}\theta_{3},\zeta
)-\Phi_{2}(x_{j};\theta_{1}\theta_{3},\zeta)\right]  -\theta_{2}\Phi_{2}%
(x_{j};\theta_{2}\theta_{3},\zeta)\Phi_{1}(x_{j};\theta_{1},\zeta)}{\theta
_{2}^{(x_{j}+1)^{-\zeta}}\Phi_{1}(x_{j};\theta_{1}\theta_{3},\zeta)-\left(
\theta_{2}\theta_{3}\right)  ^{(x_{j})^{-\zeta}}\Phi_{1}(x_{j};\theta
_{1},\zeta)},\label{A7}%
\end{align}
and%
\begin{align}
\frac{\partial L}{\partial\zeta}  & =\overset{n_{1}}{\underset{j=1}{\sum}%
}\frac{\Phi_{3}(x_{1j};\theta_{1}\theta_{3},\zeta)-\Phi_{3}(x_{1j}%
+1;\theta_{1}\theta_{3},\zeta)}{\Phi_{1}(x_{1j};\theta_{1}\theta_{3},\zeta
)}+\overset{n_{1}}{\underset{j=1}{\sum}}\frac{\Phi_{3}(x_{2j};\theta_{2}%
,\zeta)-\Phi_{3}(x_{2j}+1;\theta_{2},\zeta)}{\Phi_{1}(x_{2j};\theta_{2}%
,\zeta)}\nonumber\\
& +\overset{n_{2}}{\underset{j=1}{\sum}}\frac{\Phi_{3}(x_{1j};\theta_{1}%
,\zeta)-\Phi_{3}(x_{1j}+1;\theta_{1},\zeta)}{\Phi_{1}(x_{1j};\theta_{1}%
,\zeta)}+\overset{n_{2}}{\underset{j=1}{\sum}}\frac{\Phi_{3}(x_{2j};\theta
_{2}\theta_{3},\zeta)-\Phi_{3}(x_{2j}+1;\theta_{2}\theta_{3},\zeta)}{\Phi
_{1}(x_{2j};\theta_{2}\theta_{3},\zeta)}\nonumber\\
& \overset{n_{3}}{+\underset{j=1}{\sum}}\frac{\theta_{2}^{(x_{j}+1)^{-\zeta}%
}\left[  \Phi_{3}(x_{j};\theta_{1}\theta_{3},\zeta)-\Phi_{3}(x_{j}%
+1;\theta_{1}\theta_{3},\zeta)\right]  -\Phi_{1}(x_{j};\theta_{1}\theta
_{3},\zeta)\Phi_{3}(x_{j}+1;\theta_{2},\zeta)}{\theta_{2}^{(x_{j}+1)^{-\zeta}%
}\Phi_{1}(x_{j};\theta_{1}\theta_{3},\zeta)-\left(  \theta_{2}\theta
_{3}\right)  ^{(x_{j})^{-\zeta}}\Phi_{1}(x_{j};\theta_{1},\zeta)}\nonumber\\
& \overset{n_{3}}{+\underset{j=1}{\sum}}\frac{\left(  \theta_{2}\theta
_{3}\right)  ^{(x_{j})^{-\zeta}}\left[  \Phi_{3}(x_{j}+1;\theta_{1}%
,\zeta)-\Phi_{3}(x_{j};\theta_{1},\zeta)\right]  -\Phi_{1}(x_{j};\theta
_{1},\zeta)\Phi_{3}(x_{j};\theta_{2}\theta_{3},\zeta)}{\theta_{2}%
^{(x_{j}+1)^{-\zeta}}\Phi_{1}(x_{j};\theta_{1}\theta_{3},\zeta)-\left(
\theta_{2}\theta_{3}\right)  ^{(x_{j})^{-\zeta}}\Phi_{1}(x_{j};\theta
_{1},\zeta)},\label{A8}%
\end{align}
where $\Phi_{2}(x;\theta,\zeta)=x^{-\zeta}\theta^{x^{-\zeta}-1}$ and $\Phi
_{3}(x;\theta,\zeta)=x^{-\zeta}\theta^{x^{-\zeta}-1}\ln(x)\ln(\theta).$ The
MLEs of the parameters $\theta_{1},\theta_{2},\theta_{3}$ and $\zeta$ can be
obtained by solving the above system of four non-linear equations from
Equation (\ref{A5}) to Equation (\ref{A8}). The solution of these equations is
not easy to solve, so we need a numerical technique to get the MLEs.

\subsection{Simulation results}

In this section, we introduce some simulation results to show how the proposed
MLE performs for different sample sizes and for different parameter values.
So, we have taken two sets of parameter values: $\theta_{1}=0.8,\theta
_{2}=0.4,\theta_{3}=0.4,\zeta=0.5$ and $\theta_{1}=0.6,\theta_{2}%
=0.25,\theta_{3}=0.3,\zeta=0.9.$ The population parameters are generated using
software "Mathcad prime 3" package. The sampling distributions are obtained
for different sample sizes $n=[50,100,150,250,400]$ from $N=500$ replications.
In each case we have generated a random sample from the BDsIW$\left(
\theta_{1},\theta_{2},\theta_{3},\zeta\right)  $ with the given sample size
and the parameter values. Tables 1 and 2 obtain the average estimates (AvE)
and the mean squared errors (MSEs) of the different parameters.%

\begin{align*}
& \text{\textbf{Table 1.} The AvE and MSE values for the BDsIW(}%
0.8,0.4,0.4,0.5\text{).}\\
&
\begin{tabular}
[c]{|c|c|c|c|c|c|c|c|c|}\hline\hline
\textbf{Size} & \multicolumn{2}{|c|}{$\theta_{1}$} &
\multicolumn{2}{|c|}{$\theta_{2}$} & \multicolumn{2}{|c|}{$\theta_{3}$} &
\multicolumn{2}{|c|}{$\zeta$}\\\hline\hline
$\mathbf{n}$ & \textbf{AvE} & \textbf{MSE} & \textbf{AvE} & \textbf{MSE} &
\textbf{AvE} & \textbf{MSE} & \textbf{AvE} & \textbf{MSE}\\\hline\hline
\textbf{50} & 0.765 & 0.0311 & 0.424 & 0.0229 & 0.399 & 0.0147 & 0.515 &
0.0038\\\hline
\textbf{100} & 0.770 & 0.0307 & 0.412 & 0.0226 & 0.398 & 0.0129 & 0.504 &
0.0017\\\hline
\textbf{150} & 0.771 & 0.0303 & 0.414 & 0.0215 & 0.399 & 0.0119 & 0.497 &
0.0010\\\hline
\textbf{250} & 0.774 & 0.0299 & 0.413 & 0.0194 & 0.402 & 0.0102 & 0.503 &
0.0005\\\hline
\textbf{400} & 0.788 & 0.0284 & 0.410 & 0.0193 & 0.401 & 0.0100 & 0.499 &
0.0004\\\hline\hline
\end{tabular}
\end{align*}%
\begin{align*}
& \text{\textbf{Table 2.} The AvE and MSE values for the BDsIW(}%
0.6,0.25,0.3,0.9\text{).}\\
&
\begin{tabular}
[c]{|c|c|c|c|c|c|c|c|c|}\hline\hline
\textbf{Size} & \multicolumn{2}{|c|}{$\theta_{1}$} &
\multicolumn{2}{|c|}{$\theta_{2}$} & \multicolumn{2}{|c|}{$\theta_{3}$} &
\multicolumn{2}{|c|}{$\zeta$}\\\hline\hline
$\mathbf{n}$ & \textbf{AvE} & \textbf{MSE} & \textbf{AvE} & \textbf{MSE} &
\textbf{AvE} & \textbf{MSE} & \textbf{AvE} & \textbf{MSE}\\\hline\hline
\textbf{50} & 0.667 & 0.0340 & 0.285 & 0.0321 & 0.295 & 0.0190 & 0.878 &
0.0070\\\hline
\textbf{100} & 0.663 & 0.0326 & 0.283 & 0.0304 & 0.295 & 0.0157 & 0.882 &
0.0055\\\hline
\textbf{150} & 0.661 & 0.0311 & 0.283 & 0.0291 & 0.297 & 0.0137 & 0.884 &
0.0032\\\hline
\textbf{250} & 0.660 & 0.0284 & 0.280 & 0.0212 & 0.293 & 0.0136 & 0.887 &
0.0023\\\hline
\textbf{400} & 0.653 & 0.0202 & 0.279 & 0.0204 & 0.290 & 0.0135 & 0.890 &
0.0017\\\hline\hline
\end{tabular}
\end{align*}
Based on the simulation results, it is observed that as $n$ increases, the MSE
decreases. Moreover, the AvE and initial values are approximately equal. So;
the MLE can be used quite effectively for data analysis purposes.

\subsection{Data analysis}

In this section, we explain the experimental importance of BDsIW distribution
using two applications to real data sets. In each data, we shall compare the
fits of BDsIW distribution with some competitive models. The tested
distributions are compared using some criteria namely, the maximized
log-likelihood (\ $-L$\ ), Akaike information criterion (AIC), corrected
Akaike information criterion (CAIC), bayesian information criterion (BIC) and
Hannan-Quinn information criterion (HQIC). Further, we can use the Pearson's
chi-square goodness-of-fit test for grouped data to test the goodness of fit
of a proposed bivariate distribution. But the sample size must be sufficiently
large in order to apply this test. For this reason, we did not use this test
in the two data sets analyzed here.

\subsubsection{The first data: Football data}

This data is reported in Lee and Cha (2015), and it represents a football
match score in Italian football match (Serie A) during 1996 to 2011, between
ACF Fiorentina($X_{1}$) and Juventus($X_{2}$). This data is reported in Table
3. \ \ %

\begin{align*}
& \text{\textbf{Table 3. }The\textbf{\ }score data between ACF Fiorentina and
Juventus.}\\
&
\begin{tabular}
[c]{|c|c|c|c||c|c|c|c|}\hline\hline
\textbf{Obs.} & \textbf{Match Date} & $X_{1}$ & $X_{2}$ & \ \textbf{Obs.} &
\textbf{Match Date} & $X_{1}$ & $X_{2}$\\\hline\hline
\textbf{1} & 25th Oct. 2011 & $1$ & $2$ & \textbf{14} & 16th Feb. 2002 & $1$ &
$2$\\\hline
\textbf{2} & 17th Apr. 2011 & $0$ & $0$ & \textbf{15} & 19th Dec. 2001 & $1$ &
$1$\\\hline
\textbf{3} & 27th Nov. 2010 & $1$ & $1$ & \textbf{16} & 12th May. 2001 & $1$ &
$3$\\\hline
\textbf{4} & 06th Mar. 2010 & $1$ & $2$ & \textbf{17} & 06th Jan. 2001 & $3$ &
$3$\\\hline
\textbf{5} & 17th Oct. 2009 & $1$ & $1$ & \textbf{18} & 21st Apr. 2000 & $0$ &
$1$\\\hline
\textbf{6} & 24th Jan. 2009 & $0$ & $1$ & \textbf{19} & 18th Dec. 1999 & $1$ &
$1$\\\hline
\textbf{7} & 31st Aug. 2008 & $1$ & $1$ & \textbf{20} & 24th apr. 1999 & $1$ &
$2$\\\hline
\textbf{8} & 02nd Mar. 2008 & $3$ & $2$ & \textbf{21} & 12th Dec. 1998 & $1$ &
$0$\\\hline
\textbf{9} & 07th Oct. 2007 & $1$ & $1$ & \textbf{22} & 21st Feb. 1998 & $3$ &
$0$\\\hline
\textbf{10} & 09th Apr. 2006 & $1$ & $1$ & \textbf{23} & 04th Oct. 1997 & $1$
& $2$\\\hline
\textbf{11} & 04th Dec. 2005 & $1$ & $2$ & \textbf{24} & 22nd Feb. 1997 & $1$
& $1$\\\hline
\textbf{12} & 09th Apr. 2005 & $3$ & $3$ & \textbf{25} & 28th Sept. 1996 & $0
$ & $1$\\\hline
\textbf{13} & 10th Nov.2004 & $0$ & $1$ & \textbf{26} & 23rd Mar. 1996 & $0$ &
$1$\\\hline\hline
\end{tabular}
\end{align*}
We shall compare the fits of BDsIW distribution with some competitive models
like BDsE, BDsR, BDsW, bivariate Poisson with minimum operator (BPo$_{\min}$),
bivariate Poisson with 3-parameter (BPo-3P), independent bivariate Poisson
(IBPo), BDsIE, and BDsIR distributions. Before trying to analyze the data
using BDsIW distribution, we fit at first the marginals $X_{1}$ and $X_{2}$
separately and the $\min(X_{1},X_{2})$ on this data. The MLEs of the
parameters $\theta$ and $\zeta$ of the corresponding DsIW distribution for
$X_{1}$, $X_{2}$ and $\min(X_{1},X_{2})$ are (0.237, 2.798), (0.095, 2.601)
and (0.310, 3.103) respectively. Moreover, the $-L$ values are 30.86, 33.73
and 28.02 respectively. Figure 3 shows the estimated PMF plots for the
marginals $X_{1}$, $X_{2}$ and $\min(X_{1},X_{2})$ using this data.%

\[%
{\includegraphics[
height=1.8343in,
width=1.8343in
]%
{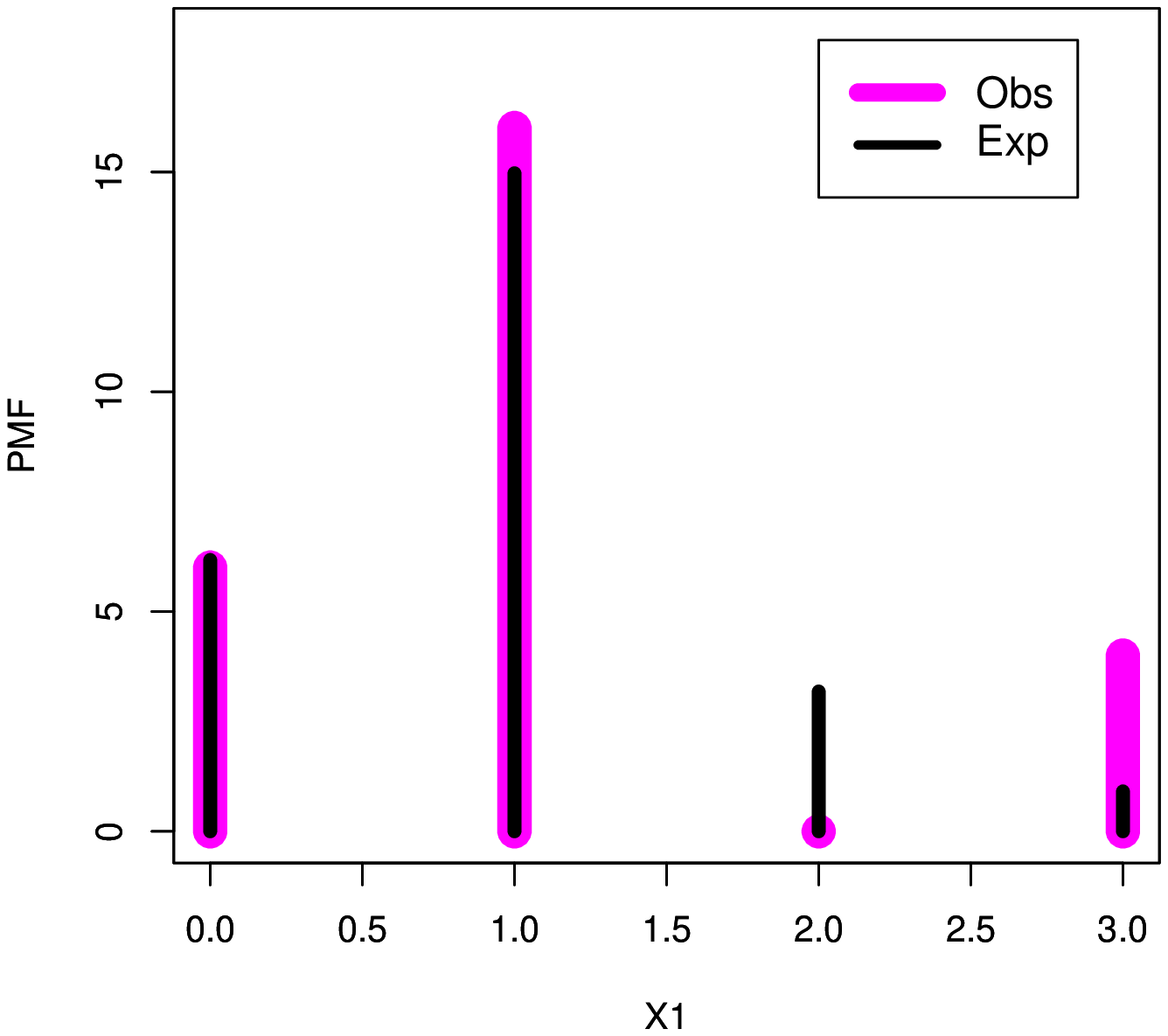}%
}
{\includegraphics[
height=1.8343in,
width=1.8343in
]%
{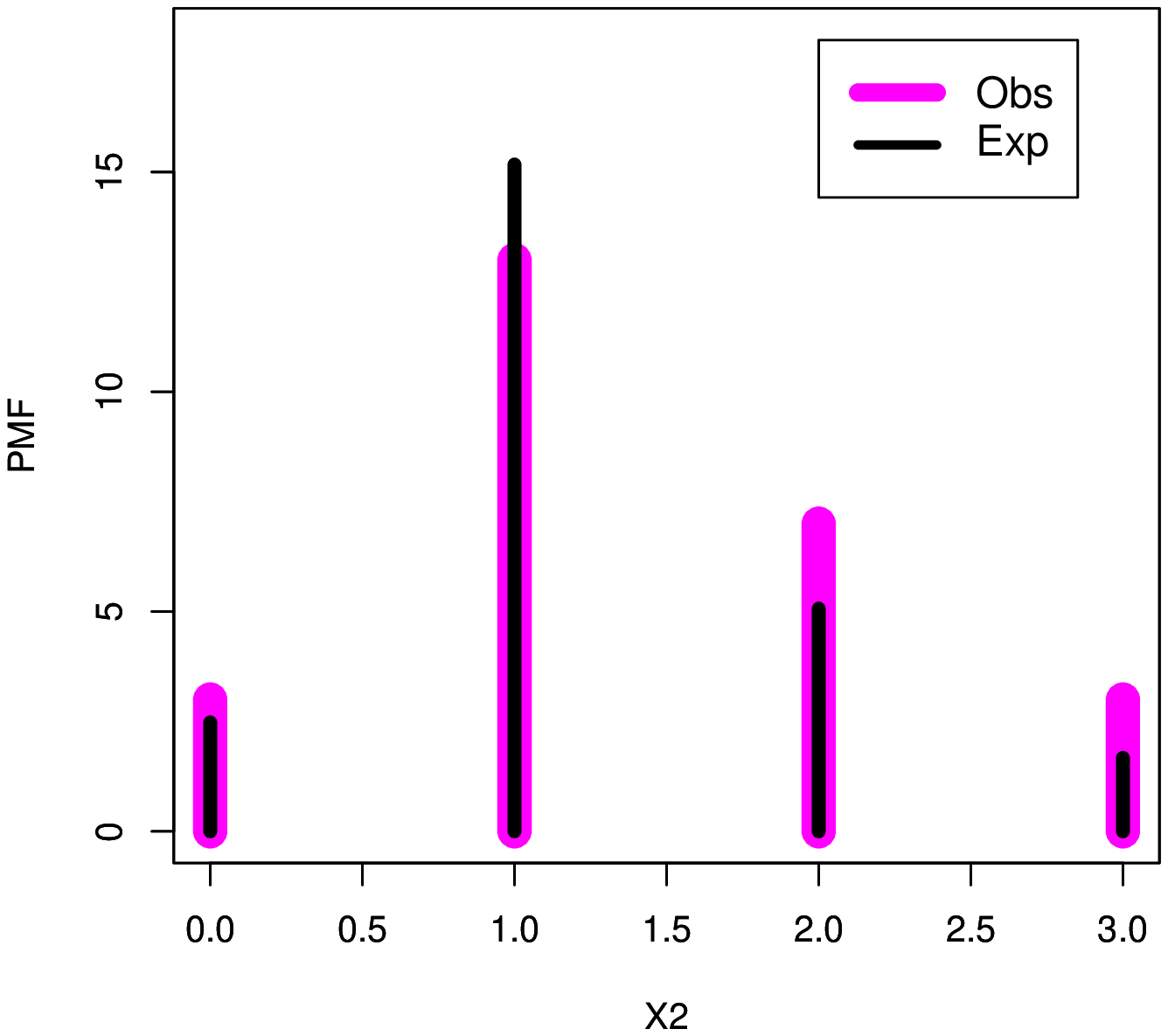}%
}
{\includegraphics[
height=1.8343in,
width=1.8343in
]%
{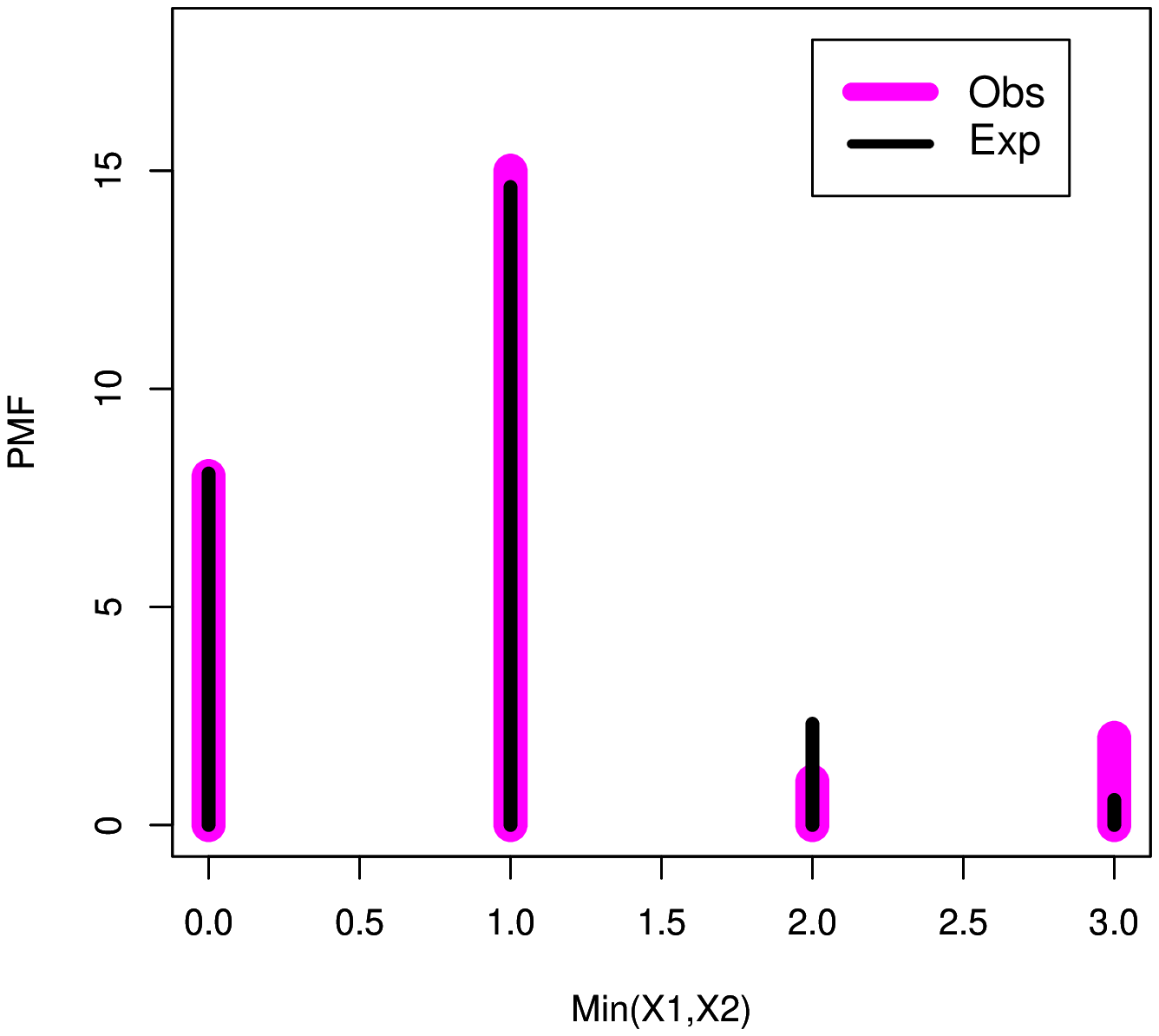}%
}
\]

\begin{center}
Figure 3. The estimated PMF for the marginals $X_{1}$, $X_{2}$ and $\min
(X_{1},X_{2})$ using football data set.
\end{center}

From Figure 3, it is clear that DsIW distribution fits the data{\small \ }for
the marginals. Now, we fit BDsIW distribution on this data. The MLEs, $-L $,
AIC, CAIC, BIC, and HQIC values for the tested bivariate models are reported
in Table 4.%

\[%
\begin{tabular}
[c]{|c|c|c|c|c|c|c|}%
\multicolumn{7}{l|}{$\text{\textbf{Table 4.} The MLEs, }-L\text{, AIC, CAIC,
BIC, and HQIC values. }$}\\\hline\hline
\textbf{Model} & \textbf{MLEs} & $\mathbf{-L}$ & \textbf{AIC} & \textbf{CAIC}
& \textbf{BIC} & \textbf{HQIC}\\\hline\hline
\multicolumn{1}{|l|}{\textbf{BDsE}} & \multicolumn{1}{|l|}{$\widehat{\theta
}_{1}=0.652,\widehat{\theta}_{2}=0.812,\widehat{\theta}_{3}=0.713$} & $75.35 $
& $156.70$ & $157.79$ & $160.47$ & $157.79$\\
\multicolumn{1}{|l|}{\textbf{BDsR}} & \multicolumn{1}{|l|}{$\widehat{\theta
}_{1}=0.790,\widehat{\theta}_{2}=0.872,\widehat{\theta}_{3}=0.905$} & $63.99 $
& $133.98$ & $135.07$ & $137.75$ & $135.07$\\
\multicolumn{1}{|l|}{\textbf{BDsW}} & \multicolumn{1}{|l|}{$\widehat{\theta
}_{1}=0.807,\widehat{\theta}_{2}=0.882,\widehat{\theta}_{3}=0.917,\widehat
{\zeta}=2.125$} & $63.89$ & $133.78$ & $134.87$ & $137.55$ & $134.87$\\
\multicolumn{1}{|l|}{\textbf{\ BPo}$_{\min}$} & \multicolumn{1}{|l|}{$\widehat
{\theta}_{1}=1.36,\widehat{\theta}_{2}=2.10,\widehat{\theta}_{3}=2.27$} &
$64.22$ & $134.44$ & $135.53$ & $138.21$ & $135.53$\\
\multicolumn{1}{|l|}{\textbf{BPo-3P}} & \multicolumn{1}{|l|}{$\widehat{\alpha
}_{1}=1.08,\widehat{\alpha}_{2}=1.38,\widehat{\alpha}_{3}=0.70$} & $64.92$ &
$135.83$ & $136.93$ & $139.61$ & $136.93$\\
\multicolumn{1}{|l|}{\textbf{IBPo}} & \multicolumn{1}{|l|}{$\widehat{\lambda
}_{1}=1.08,\widehat{\lambda}_{2}=1.38$} & $67.60$ & $139.21$ & $139.72$ &
$141.72$ & $139.92$\\
\multicolumn{1}{|l|}{\textbf{BDsIE}} & \multicolumn{1}{|l|}{$\widehat{\theta
}_{1}=0.669,\widehat{\theta}_{2}=0.388,\widehat{\theta}_{3}=0.514$} & $78.54$
& $163.07$ & $163.99$ & $167.28$ & $164.42$\\
\multicolumn{1}{|l|}{\textbf{BDsIR}} & \multicolumn{1}{|l|}{$\widehat{\theta
}_{1}=0.493,\widehat{\theta}_{2}=0.212,\widehat{\theta}_{3}=0.561$} & $64.10$
& $134.21$ & $135.29$ & $137.98$ & $135.29$\\
\multicolumn{1}{|l|}{\textbf{BDsIW}} & \multicolumn{1}{|l|}{$\widehat{\theta
}_{1}=0.420,\widehat{\theta}_{2}=0.141,\widehat{\theta}_{3}=0.587,\widehat
{\zeta}=2.738$} & $\mathbf{61.96}$ & $\mathbf{131.82}$ & $\mathbf{133.82}$ &
$\mathbf{136.95}$ & $\mathbf{133.37}$\\\hline\hline
\end{tabular}
\]

\ From Table 4, it is clear that BDsIW distribution provides a better fit than
the other tested distributions, because it has the smallest values among $-L$,
AIC, CAIC, BIC and HQIC. Since, BDsIE and BDsIR\ distributions are special
cases from BDsIW distribution. Hence, we want to perform the following two tests:

\begin{enumerate}
\item[ \ \ \ Test 1:] $H_{01}:$ $\zeta=1$ (BDsIE) against $H_{11}:$ $\zeta
\neq1$ (BDsIW).

\item[ \ \ \ Test 2:] $H_{02}:$ $\zeta=2$ (BDsIR) against $H_{12}:$ $\zeta
\neq2$ (BDsIW).
\end{enumerate}

The likelihood ratio test statistics ($\Lambda$), d.f and p-values for BDsIE
and BDsIR distributions are given in Table 5.%

\[%
\begin{tabular}
[c]{ccccc}%
\multicolumn{5}{l}{$\text{\textbf{Table 5.} The }\Lambda\text{, d.f and
p-values.}$}\\\hline\hline
{\small \ \ \textbf{\ Model \ }\ } & \ \ \ {\small \ }$\ \ \ \ H_{\circ}%
${\small \ \ \ \ \ \ \ \ \ \ \ } & {\small \ }$\ \ \ \ \ \ \ \ \ \ \Lambda
${\small \ \ \ \ \ \ \ \ \ \ \ } & \textbf{\ \ \ \ \ \ d.f. \ \ \ \ \ \ } &
\textbf{\ \ \ \ \ \ \ \ p-values \ \ \ \ \ \ \ }\\\hline\hline
\textbf{BDsIE} & $\zeta=1$ & $33.152$ & $1$ & $<0.01$\\
\textbf{BDsIR} & $\zeta=2$ & $4.288$ & $1$ & $0.0384$\\\hline\hline
\end{tabular}
\]
We can conclude that $H_{01}$ and $H_{02}$ are rejected with 5\% level of
significance. Hence, BDsIE and BDsIR distributions cannot be used for this
data set. So, we prefer BDsIW distribution for analyzing this data. Figure 4
shows the estimated joint PMF for BDsIW, BDsIE and BDsIR distributions using
this data, which support the results of Table 5.

\begin{center}%
\[%
{\includegraphics[
height=2.0349in,
width=2.0349in
]%
{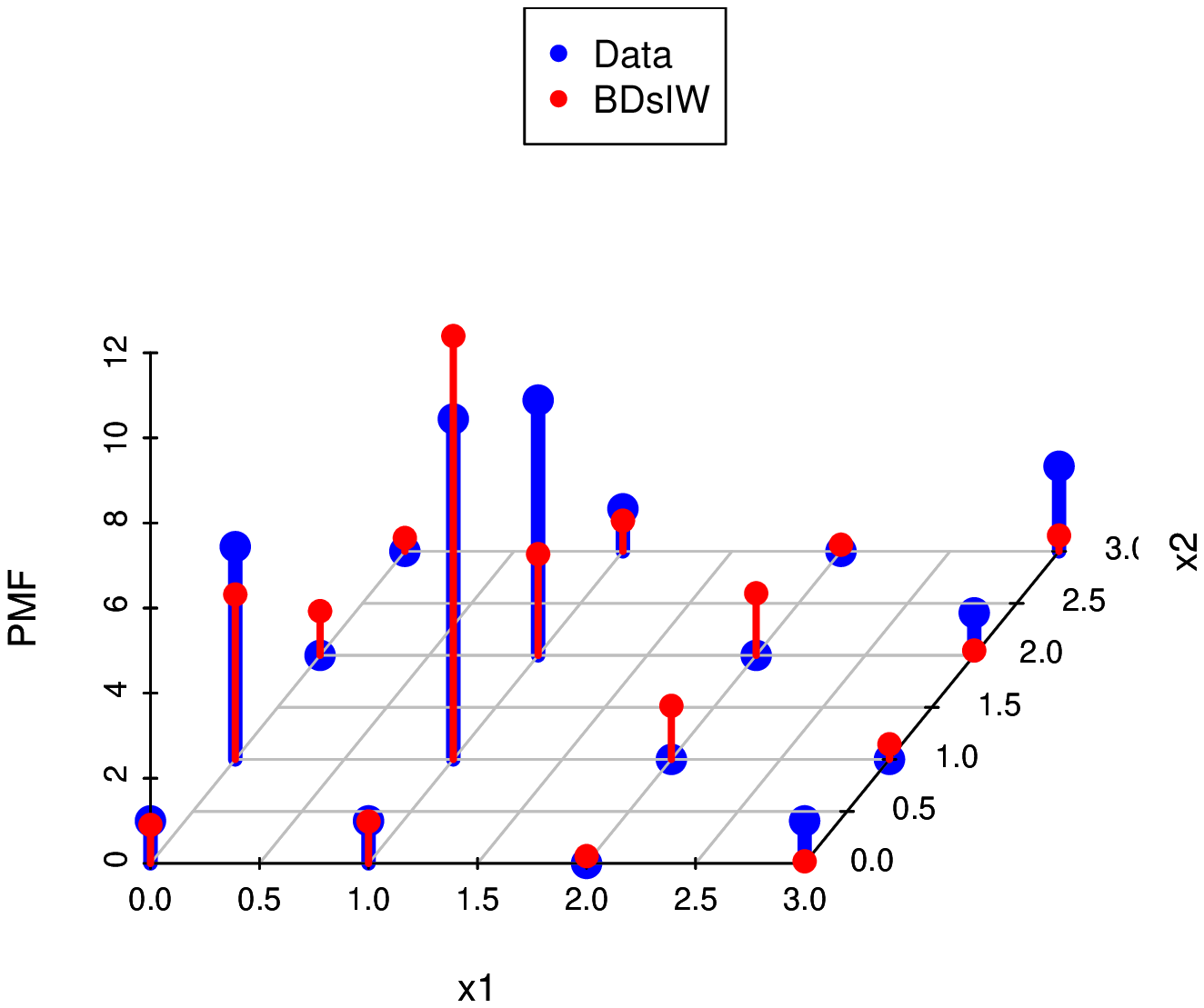}%
}
{\includegraphics[
height=2.0349in,
width=2.0349in
]%
{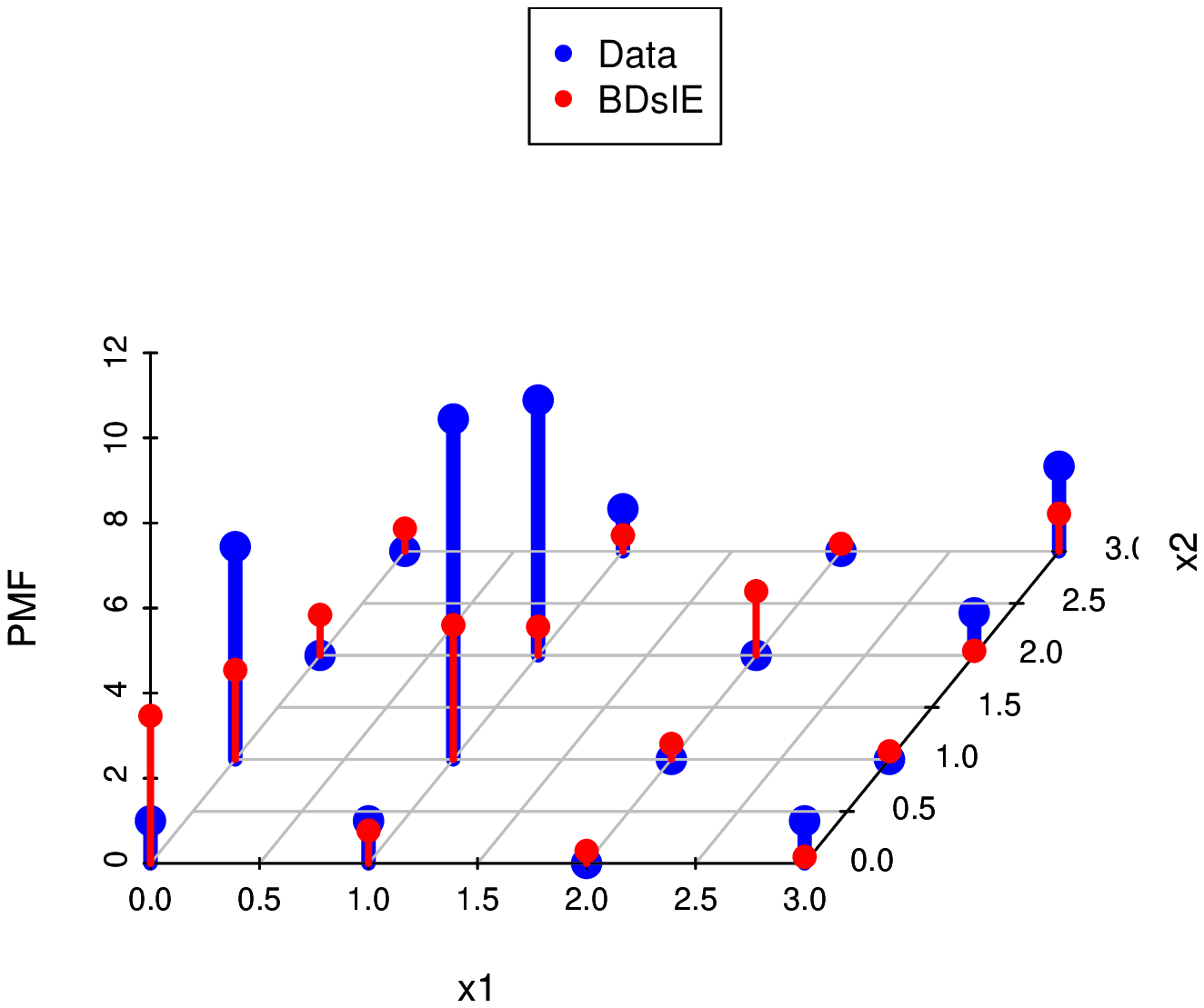}%
}
{\includegraphics[
height=2.0349in,
width=2.0349in
]%
{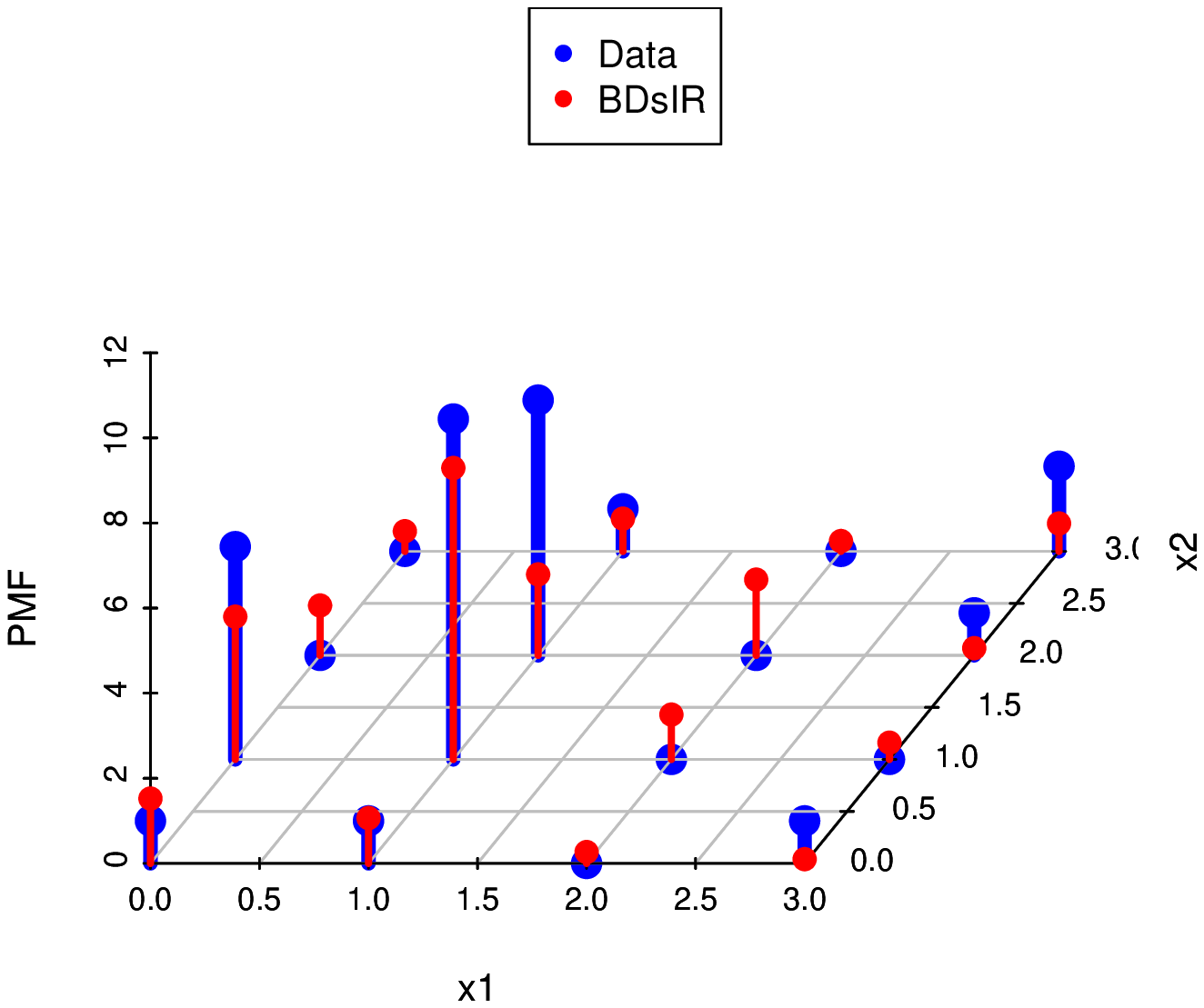}%
}
\]
Figure 4. The estimated joint PMF for BDsIW, BDsIE and BDsIR distributions
using football data set.
\end{center}

\subsubsection{The second data: Nasal drainage severity score}

This data is reported in Davis (2002), and it represents the efficacy of steam
inhalation in the treatment of common cold symptoms (0 = no symptoms; 1 = mild
symptoms; 2 = moderate symptoms; 3 = severe symptoms). This data is presented
in Table 6.%
\begin{align*}
& \text{\textbf{Table 6. }Nasal drainage severity score.}\\
&
\begin{tabular}
[c]{|c|c|c||c|c|c|}\hline\hline
\textbf{Obs.} & \textbf{Day 1 }$\mathbf{(}X_{1}\mathbf{)}$ & \textbf{Day 2}
$({\small X}_{2})$ & \textbf{Obs.} & \textbf{Day 1} $({\small X}_{1})$ &
\textbf{Day 2} $({\small X}_{2})$\\\hline\hline
\textbf{1} & $1$ & $1$ & \textbf{16} & $2$ & $1$\\\hline
\textbf{2} & $0$ & $0$ & \textbf{17} & $1$ & $1$\\\hline
\textbf{3} & $1$ & $1$ & \textbf{18} & $2$ & $2$\\\hline
\textbf{4} & $1$ & $1$ & \textbf{19} & $3$ & $1$\\\hline
\textbf{5} & $0$ & $2$ & \textbf{20} & $1$ & $1$\\\hline
\textbf{6} & $2$ & $0$ & \textbf{21} & $2$ & $1$\\\hline
\textbf{7} & $2$ & $2$ & \textbf{22} & $2$ & $2$\\\hline
\textbf{8} & $1$ & $1$ & \textbf{23} & $1$ & $1$\\\hline
\textbf{9} & $3$ & $2$ & \textbf{24} & $2$ & $2$\\\hline
\textbf{10} & $2$ & $2$ & \textbf{25} & $2$ & $0$\\\hline
\textbf{11} & $1$ & $0$ & \textbf{26} & $1$ & $1$\\\hline
\textbf{12} & $2$ & $3$ & \textbf{27} & $0$ & $1$\\\hline
\textbf{13} & $1$ & $3$ & \textbf{28} & $1$ & $1$\\\hline
\textbf{14} & $2$ & $1$ & \textbf{29} & $1$ & $1$\\\hline
\textbf{15} & $2$ & $3$ & \textbf{30} & $3$ & $3$\\\hline\hline
\end{tabular}
\end{align*}
We shall compare the fits of BDsIW distribution with some competitive models
like bivariate Poisson with 4-parameter (BPo-4P), IBPo, BDsE, BDsIE and BDsIR
distributions. We fit at first the marginals $X_{1}$ and $X_{2}$ separately
and the $\min(X_{1},X_{2})$ on this data. The MLEs of the parameters $\theta$
and $\zeta$ of the corresponding DsIW distribution for $X_{1}$, $X_{2}$ and
$\min(X_{1},X_{2})$ are (0.065, 2.505), (0.115, 2.524) and (0.181, 2.699)
respectively. Moreover, the $-L$ values are 40.99, 39.83 and 36.68
respectively. Figure 5 shows the estimated PMF plots for the marginals $X_{1}%
$, $X_{2}$ and $\min(X_{1},X_{2})$ using this data.

\begin{center}%
\[%
{\includegraphics[
height=1.8343in,
width=1.8343in
]%
{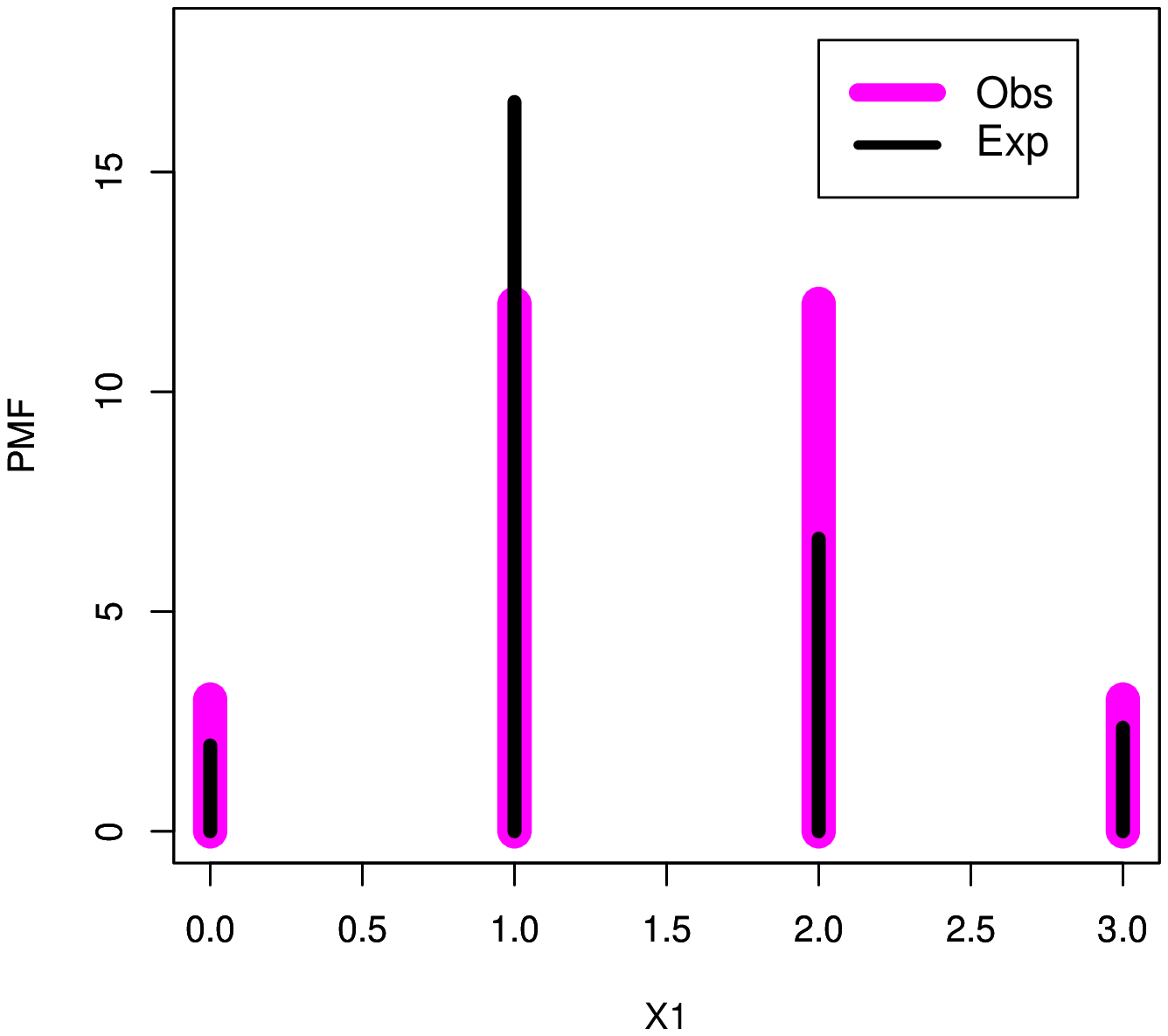}%
}
{\includegraphics[
height=1.8343in,
width=1.8343in
]%
{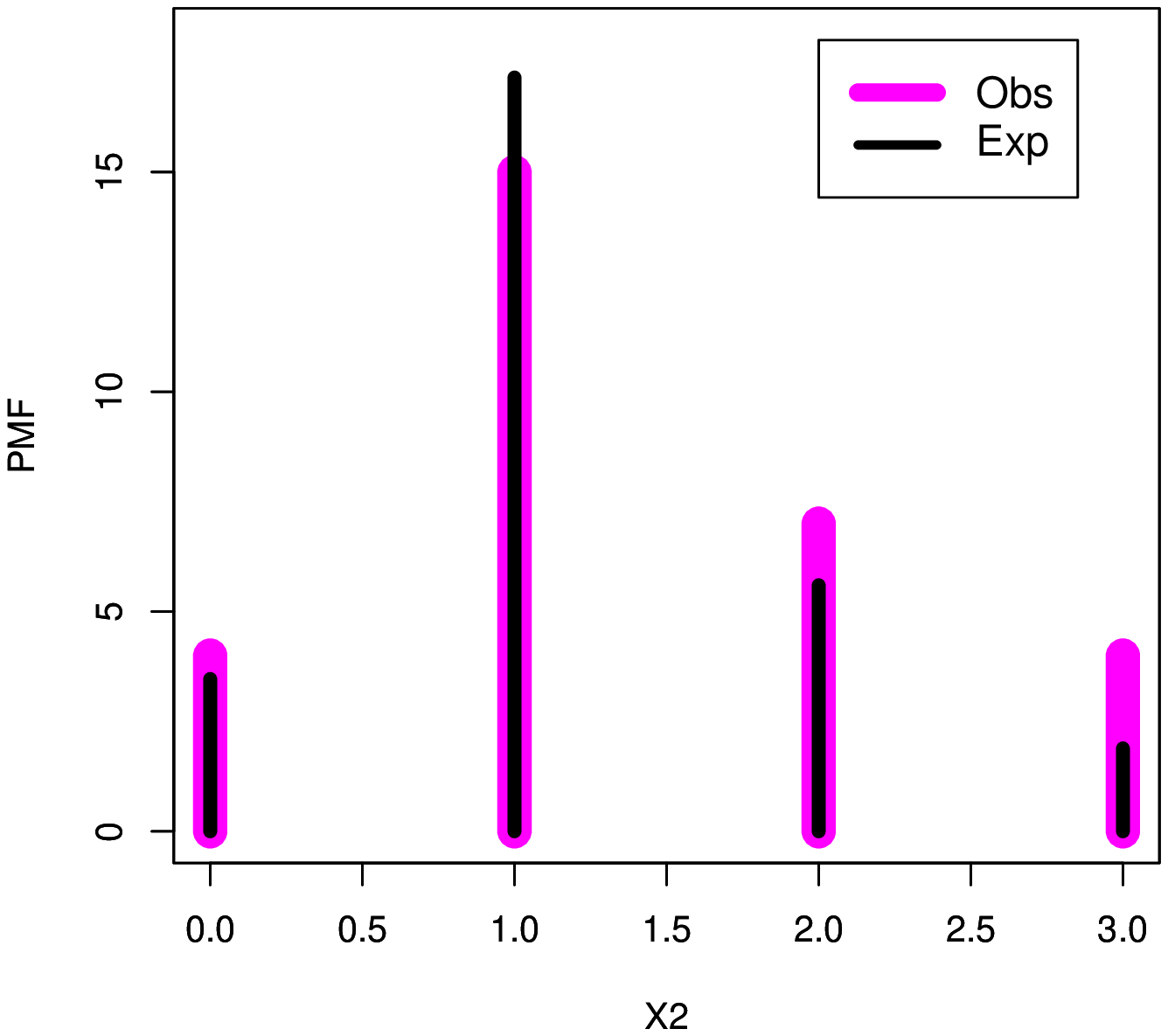}%
}
{\includegraphics[
height=1.8343in,
width=1.8343in
]%
{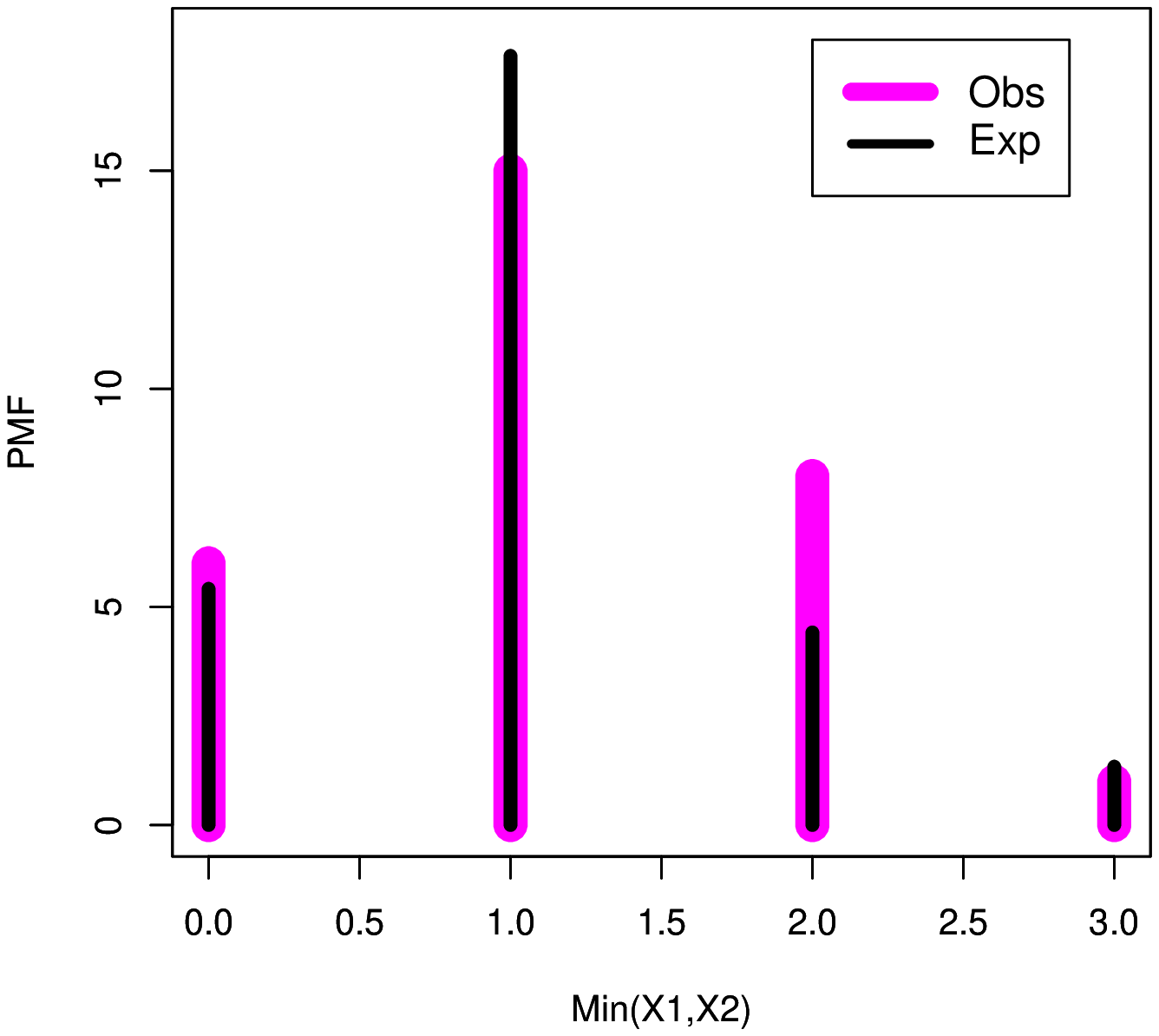}%
}
\]
Figure 5. The estimated PMF for the marginals $X_{1}$, $X_{2}$ and $\min
(X_{1},X_{2})$ using nasal drainage severity score.
\end{center}

From Figure 5, it is clear that DsIW distribution fits the data{\small \ }for
the marginals. Now, we fit BDsIW distribution on this data. The MLEs, $-L $,
AIC, CAIC, BIC, and HQIC values for the tested bivariate models are reported
in Table 7.
\begin{align*}
& \text{\textbf{Table 7.} The MLEs, }-L\text{, AIC, CAIC, BIC, and HQIC
values.}\\
&
\begin{tabular}
[c]{|c|c|c|c|c|c|c|}\hline\hline
\textbf{Model} & \textbf{MLEs} & $\mathbf{-L}$ & \textbf{AIC} & \textbf{CAIC}
& \textbf{BIC} & \textbf{HQIC}\\\hline\hline
\multicolumn{1}{|l|}{\textbf{BPo-4P}} & \multicolumn{1}{|l|}{$\widehat
{\lambda}_{1}=0.262,\widehat{\alpha}_{1}=0.165,\widehat{\lambda}%
_{2}=0.405,\widehat{\alpha}_{2}=2.97$} & $77.66$ & $163.33$ & $164.93$ &
$168.93$ & $164.66$\\
\multicolumn{1}{|l|}{\textbf{IBPo}} & \multicolumn{1}{|l|}{$\widehat{\lambda
}_{1}=1.499,\widehat{\lambda}_{2}=1.367$} & $92.48$ & $190.96$ & $191.88$ &
$195.16$ & $192.30$\\
\multicolumn{1}{|l|}{\textbf{BDsE}} & \multicolumn{1}{|l|}{$\widehat{\theta
}_{1}=0.846,\widehat{\theta}_{2}=0.792,\widehat{\theta}_{3}=0.693$} & $88.00 $
& $182$ & $182.92$ & $186.20$ & $183.34$\\
\multicolumn{1}{|l|}{\textbf{BDsIE}} & \multicolumn{1}{|l|}{$\widehat{\theta
}_{1}=0.501,\widehat{\theta}_{2}=0.622,\widehat{\theta}_{3}=0.383$} & $92.48$
& $190.96$ & $191.88$ & $195.16$ & $192.30$\\
\multicolumn{1}{|l|}{\textbf{BDsIR}} & \multicolumn{1}{|l|}{$\widehat{\theta
}_{1}=0.262,\widehat{\theta}_{2}=0.405,\widehat{\theta}_{3}=0.363$} & $78.66$
& $163.32$ & $164.24$ & $167.52$ & $164.66$\\
\multicolumn{1}{|l|}{\textbf{BDsIW}} & \multicolumn{1}{|l|}{$\widehat{\theta
}_{1}=0.192,\widehat{\theta}_{2}=0.337,\widehat{\theta}_{3}=0.360,\widehat
{\zeta}=2.453$} & $\mathbf{76.51}$ & $\mathbf{161.02}$ & $\mathbf{162.62}$ &
$\mathbf{166.62}$ & $\mathbf{162.81}$\\\hline\hline
\end{tabular}
\end{align*}
From Table 7, it is clear that BDsIW distribution provides a better fit than
the other tested distributions. Table 8 shows the $\Lambda$ and p-values for
BDsIE and BDsIR distributions using nasal drainage severity score data set.%

\[%
\begin{tabular}
[c]{ccccc}%
\multicolumn{5}{l}{$\text{\textbf{Table 8.} The }\Lambda\text{, d.f and
p-values.}$}\\\hline\hline
{\small \ \ \textbf{\ Model \ }\ } & \ \ \ {\small \ }$\ \ \ \ H_{\circ}%
${\small \ \ \ \ \ \ \ \ \ \ \ } & {\small \ }$\ \ \ \ \ \ \ \ \ \ \Lambda
${\small \ \ \ \ \ \ \ \ \ \ \ } & \textbf{\ \ \ \ \ \ d.f. \ \ \ \ \ \ } &
\textbf{\ \ \ \ \ \ \ \ p-values \ \ \ \ \ \ \ }\\\hline\hline
\textbf{BDIE} & $\zeta=1$ & $31.94$ & $1$ & $<0.01$\\
\textbf{BDIR} & $\zeta=2$ & $4.3$ & $1$ & $0.0381$\\\hline\hline
\end{tabular}
\]
From Table 8, we can conclude that $H_{01}$ and $H_{02}$ are rejected with 5\%
level of significance. So, we prefer BDsIW distribution for analyzing this
data. Figure 6 shows the estimated joint PMF for BDsIW, BDsIE and BDsIR
distributions using this data, which support the results of Table 8.

\begin{center}%
\[%
{\includegraphics[
height=2.0349in,
width=2.0349in
]%
{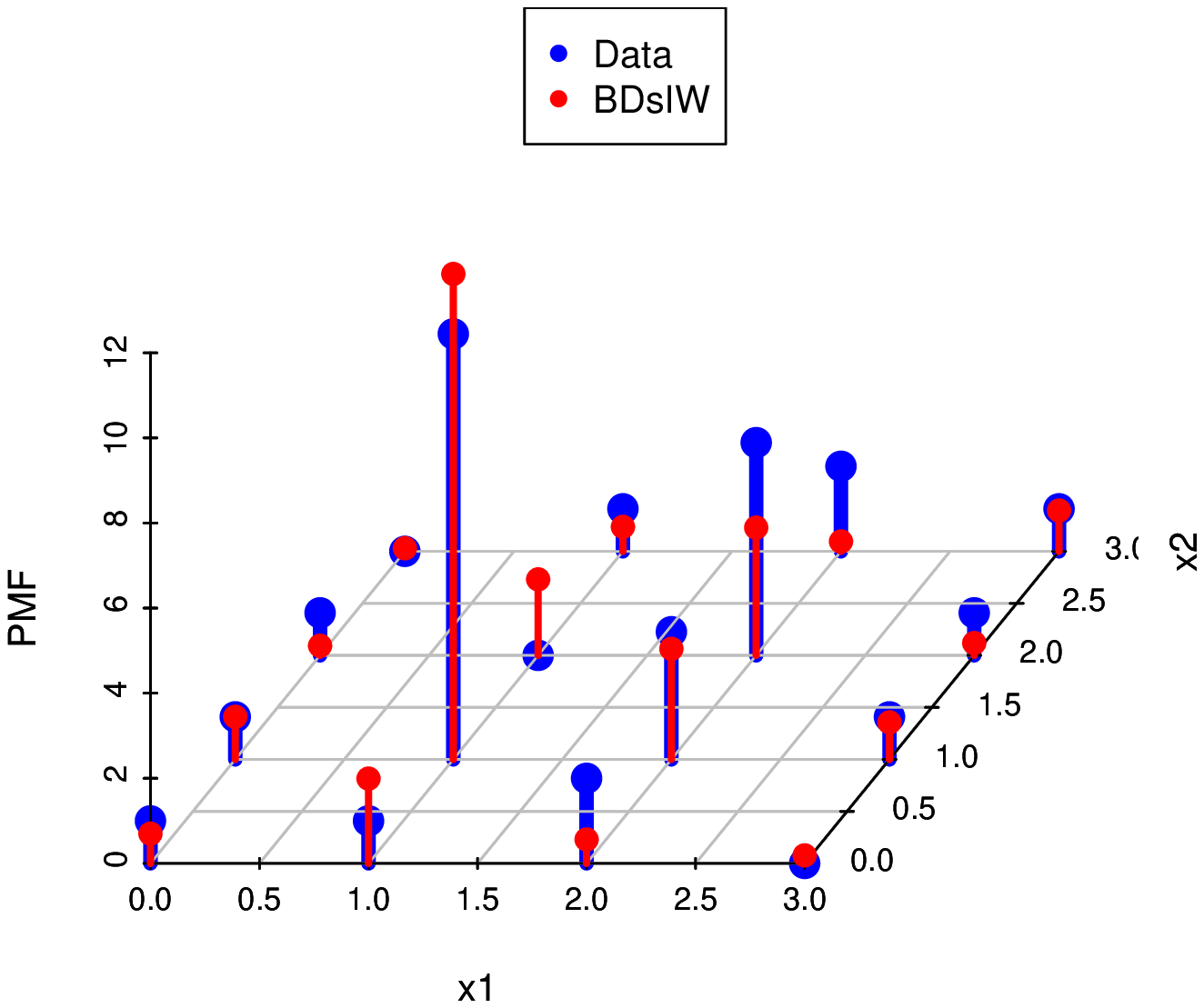}%
}
{\includegraphics[
height=2.0349in,
width=2.0349in
]%
{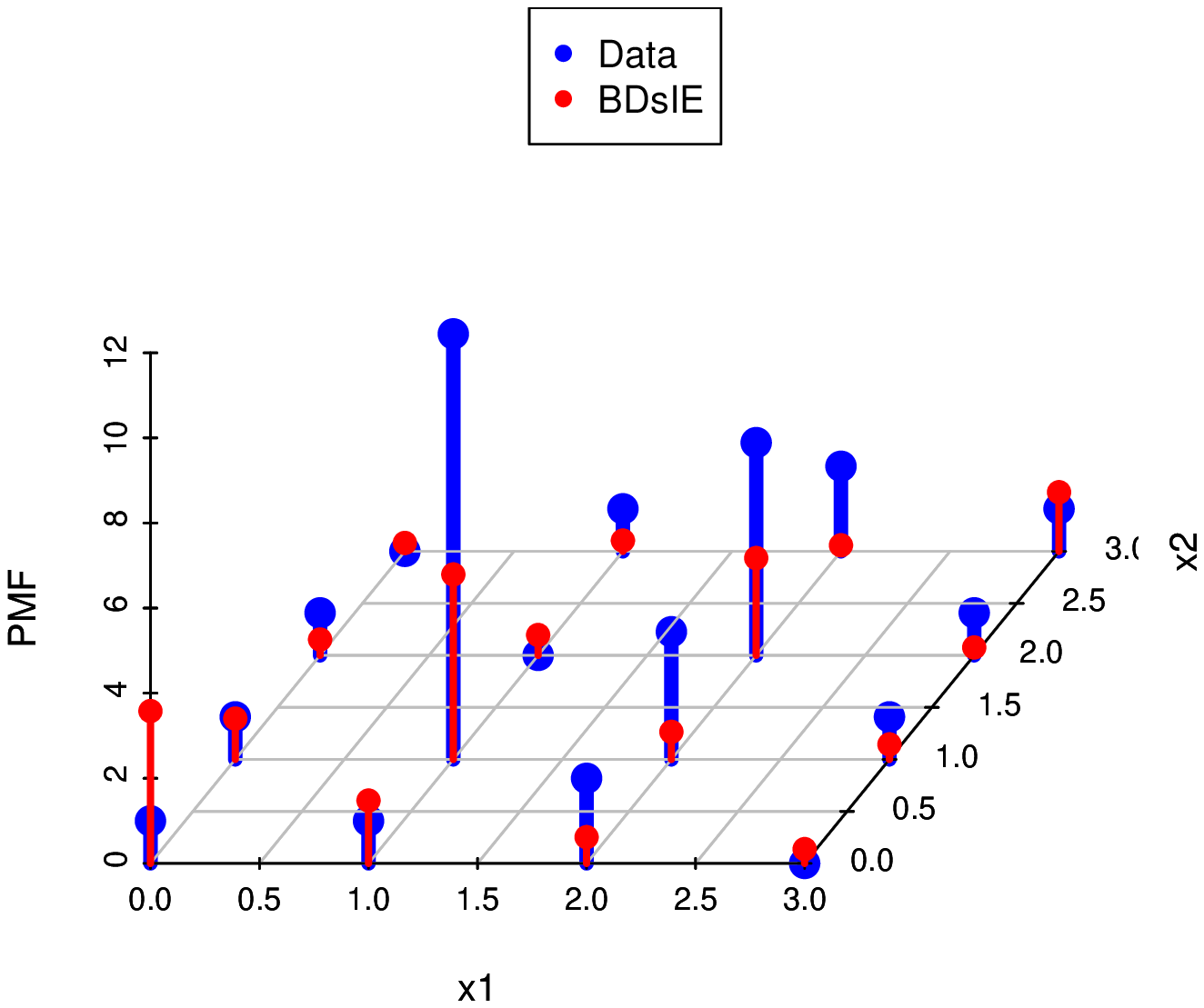}%
}
{\includegraphics[
height=2.0349in,
width=2.0349in
]%
{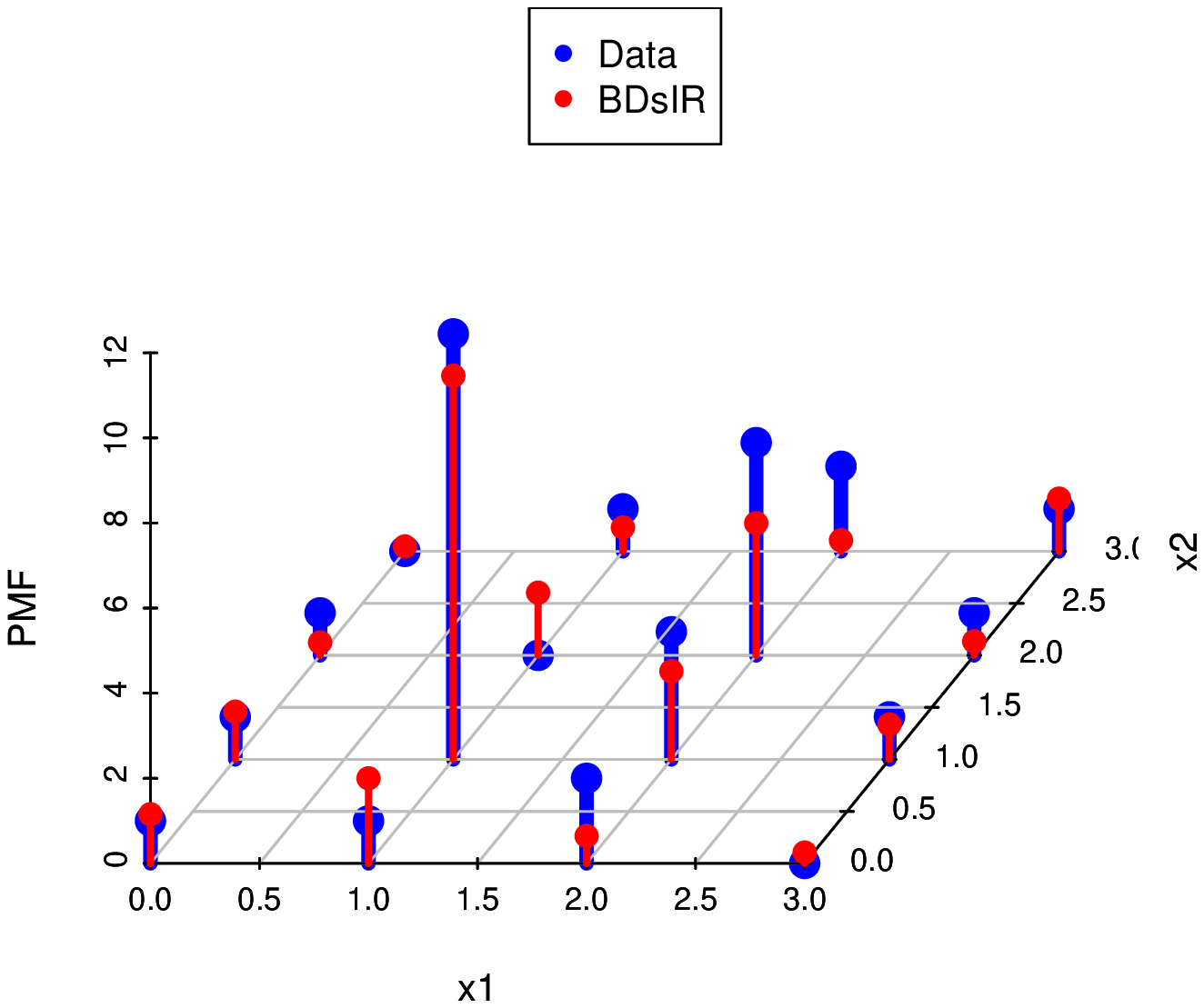}%
}
\]
Figure 6. The estimated joint PMF for BDsIW, BDsIE and BDsIR distributions
using nasal drainage severity score.
\end{center}

\section{Conclusions}

In this paper, we presented a flexible bivariate discrete distribution called
BDsIW distribution. The proposed model has the marginals, which are DsIW
distributions. The joint CDF and joint PMF have simple forms; therefore, this
new discrete model can be easily used in practice for modelling bivariate
discrete data. Some statistical and mathematical properties of the proposed
discrete model are studied. Moreover, the simulation results indicated that
the MLE works quite satisfactorily and it can be used to estimate the model
parameters. Also, we analyzed two real data sets and showed through
goodness-of-fit tests that BDsIW distribution works quite well in practice in
different fields.

\end{document}